\documentclass{amsart}
\usepackage{amssymb}
\usepackage{amscd}
\usepackage{verbatim}
\begin{document}

\newcommand{\cal}{\mathcal}
\newcommand{\s}{\bar{s} }
\newcommand{\y}{\bar{y}}
\newcommand{\bsigma}{\bar{\sigma}}
\newcommand{\bareta}{\bar{\eta}}
\newcommand{\bDelta}{\bar\Delta}

\newcommand{\R}{{\mathbb R}}
\newcommand{\C}{{\mathbb C}}
\newcommand{\Q}{{\mathbb Q}}
\newcommand{\Z}{{\mathbb Z}}
\newcommand{\half}{{\frac{1}{2}}}

\newcommand{\bz}{\bar{z}}
\newcommand{\bO}{\partial \Omega}

\newtheorem{theo}{Theorem}[section]
\newtheorem{prop}[theo]{Proposition}
\newtheorem{defn}[theo]{Definition}
\newtheorem{lem}[theo]{Lemma}
\newtheorem{cor}[theo]{Corollary}

\title[Spectral determination of analytic bi-axisymmetric plane domains ]{Spectral determination of analytic bi-axisymmetric plane domains}

\author{Steve Zelditch}
\address{Department of Mathematics, Johns Hopkins University, 
Baltimore, MD
21218, USA}
\email{zel\@math.jhu.edu }

\thanks{Research 
partially supported by  NSF grant \#DMS-9703775}

\date{\today}

\begin{abstract}  Let ${\cal D}$ denote the class of bounded real analytic
plane domains with the symmetry of an ellipse.  Under generic conditions, 
we prove that if $\Omega_1,
\Omega_2 \in {\cal D}$ and if the Dirichlet spectra coincide, $Spec(\Omega_1) =
Spec(\Omega_2)$, then $\Omega_1 = \Omega_2$ up to rigid motion. \end{abstract}

\maketitle
\centerline{\today}

\section{Introduction}

In this paper we give a positive solution to the inverse spectral problem for the
 class ${\cal D}$ of analytic axi-symmetric plane domains $\Omega$ satisfying:
$$     \left\{ \begin{array}{lll} & \bullet & \Omega\;
 \rm{ is \;real \;
analytic} \\  & \bullet &\Omega \; \rm{is}\;\Z_2 \times \Z_2-\;
\rm{symmetric} \\&  \bullet & {\rm at \;least\; one \;axis\; of }\; \Omega \;\rm{is\;
a\;non-degenerate\;bouncing\;ball\;orbit} \; \gamma\\ 
& \bullet & L_{\gamma}\;\rm{has\;multiplicity\;one\;in}\; Lsp(\Omega)\;. \end{array} \right.$$
Let Spec$(\Omega)$ denote the spectrum of the Laplacian $\Delta_{\Omega}$  of the
domain $\Omega$ with Dirichlet boundary conditions. 
Our main result is the following theorem (announced in \cite{Z.5}):
\begin{theo} \label{TH}
Spec$: {\cal D} \mapsto \R_+^{{\bf N}}$ is 1-1. \end{theo} 

Let us clarify the assumptions.  The symmetry assumption is that 
 there is an up/down reflection symmetry across a horizontal axis and a left/right reflection symmetry
across a   vertical axis.   Both axes intersect the boundary at right angles,  hence are projections to
$\Omega$ of `bouncing ball orbits'  of the the billiard flow
$G^t$ on $T^*\Omega$, with the usual law of reflection at the boundary. Associated
to any periodic reflecting ray  $\beta$ of $G^t$ is its Poincare map ${\cal P}_{\beta}$, defined as the first return map on a local transversal to $\beta$ in $S^* \Omega,$  and its  linear Poincare map $P_{\beta} = d_{\beta(0)}{\cal P}_{\beta} $.   A periodic orbit $\beta$ is said to be non-degenerate if no eigenvalue of $P_{\beta}$ is a root of unity. Thus, it is
non-degenerate elliptic if the eigenvalues
of $P_{\beta}$ have the form
$\{e^{ \pm i \alpha}\}$ with $\alpha/\pi \notin \Q$, or non-degenerate hyperbolic if they have the form
$\{e^{\pm \lambda}\}$ with $\lambda \not= 0.$  Our non-degeneracy assumption is that at least one of the axes,
which we will denote by $\overline{AB}$, 
is the projection to $\Omega$ of a non-degenerate  bouncing ball orbit $\gamma$.
Also,   $Lsp(\Omega)$ denotes the length spectrum of $\Omega,$ i.e. the set of lengths of periodic billiard trajectories
(including the boundary).   Our length spectrum assumption is  that $\gamma$ is of multiplicity
one in $Lsp(\Omega)$, i.e. that it is the unique trajectory of its length $ L_{\gamma}.$ With no loss of generality
we will assume $\overline{AB}$ is the vertical axis, and  will denote its length by  $L$;
thus $L_{\gamma} = 2 L.$

The proof of Theorem (\ref{TH}) is based on the method of normal forms, which was
introduced into inverse spectral theory by Colin de Verdiere \cite{CV} and  by Guillemin \cite{G}. The basic idea of Guillemin \cite{G} is to construct a quantum analogue of the Birkhoff normal form of $\Delta_{\Omega}$ around each closed geodesic and to prove that the coefficients
$B_{\gamma k}$ of the normal form are spectral invariants.  The latter is accomplished by relating these coefficients  $B_{\gamma k}$ to the wave invariants $a_{\gamma k}$ of $\Delta_{\Omega}$ at $\gamma$, i.e. the
coefficients of singularity expansion of the trace of the wave group $U(t)$ at the time
$t = L_{\gamma}$. In 
\cite{G}, Guillemin constructed the normal form and proved that it is a spectral invariant in the case of non-degenerate elliptic closed geodesics on boundaryless manifolds.  
A somewhat different construction of  the normal form and a new proof that the
coefficients are spectral invariants was given in
\cite{Z.1} \cite{Z.2} for general non-degenerate closed geodesics.  The inverse
spectral problem is then reduced to  determining the metric from the normal form. 
 
The latter inverse problem  remains difficult since the different closed geodesics do not easily `communicate' with each other and since the $B_{\gamma k}$'s for a fixed $\gamma$
(and its iterates) do not appear to give enough information to determine the metric, even
locally.  Therefore it is natural to  consider the problem first for classes
of real analytic metrics with one functional degree of freedom, where it is plausible
that the normal form coefficients at just one closed geodesic should determine the metric.  This
motivated our study in 
\cite{Z.3}  of real analytic `simple' surfaces of revolution.  We
proved there  that a simple analytic surface of revolution is determined by the normal form 
of its $\Delta_{\Omega}$ along the `meridian torus' of closed
geodesics.  For further background on the method of normal forms we refer to the expository artice \cite{Z.4}.

   In this paper, we extend the method of normal forms to the case of bounded analytic  plane
domains with a non-degenerate bouncing ball orbit $\gamma$.  We do not use the symmetry assumption on the domain
in the construction of the normal form at $\gamma$, but only in the last step of deducing the domain from the normal form.
Since the construction of the normal form and the various technical pitfalls may make  the proof
difficult to follow, let us give a brief summary here of the main ideas. See also
\cite{Z.5} for a somewhat less technical outline.  

 As mentioned above, the
main  idea is to   introduce a notion 
of Birkhoff normal
form $F(|D|, \hat{I})^2$  for $\Delta$ at a bouncing ball orbit $\gamma$. The normal form lives on the model
space  $\Omega_o : = [0, L]_{s} \times \R_{y}$, which carries a natural abelian algebra of pseudodifferential operators ${\cal A} = \langle |D_s|, \hat{I} \rangle$, where:
$|D|= |D_s| = \sqrt{- D_s^* D_s}$ with  Dirichlet
boundary conditions on $[0, L]$ and where $\hat{I}$ is a quantum action operator. Here, $D_x = \frac{1}{i} \frac{\partial}{\partial x}$. In the elliptic case, 
$\hat{I} = \half ( |D|^{-1} D_y^2 + |D| y^2)$ is the `transverse' homogeneous harmonic oscillator, while in the
hyperbolic case $\hat{I} = \half (|D|^{-1} D_y^2 - |D| y^2)$. By definition, the normal form of $\sqrt{\Delta}$   is  a first order polyhomogeneous symbol
\begin{equation} \label{NF} F(|D|, \hat{I}) \sim |D| + \frac{\alpha}{L} \hat{I} + \frac{p_1(\hat{I})}{|D|} + \frac{p_2(\hat{I})}{|D|^2} + \dots \end{equation}
in $|D|, \hat{I}$
with $p_j$ a polynomial of degree $j + 1.$ The coefficients $B_{\gamma k}$ mentioned above
are the coefficients of these polynomials (cf. \cite{G} \cite{Z.1} \cite{Z.4}). 

 Heuristically, $\sqrt{\Delta}$ should be microlocally
conjugate to $F(|D|, \hat{I})$ near corresponding bouncing ball orbits.  This should imply that  $\sqrt{\Delta}$
and $F(|D|, \hat{I})$ have the same wave invariants at iterates of their corresponding bouncing ball orbits, and consequently that the formal function $F$ is a spectral invariant.  
As may be anticipated the boundary gives
rise to many complications and we cannot quite implement this outline. In fact, we  work as much as possible in the open space containing
the bounded domain. The normal form is actually used to define a parametrix for the Dirichlet wave kernel near $\gamma$.

\subsection{Outline of the proof }

The main contribution of this paper is the construction of a normal form $F(|D|, \hat{I})$  of  $\Delta$
around a bouncing ball orbit $\gamma$, and the proof that  the normal form is a spectral invariant.  Here are the main steps.   

\subsubsection{Straightening the domain}

 To conjugate to the normal form,
we first use  a special  map $\Phi$ introduced by Lazutkin in \cite{L} 
to `straighten the domain' near $\overline{AB}$, i.e.
to carry an open neighborhood of $\gamma$ in $\Omega$ to an open neighborhood
of $[0, L] \times \{0\}$ in $\Omega_o$. In fact, $\Phi$ will be defined in a
neighborhood in $\R^2$ of $\overline{AB}.$  In the elliptic case, Lazutkin's map additionally puts the
metric into the  normal form $ds^2 + b(s, y) [y^2 ds^2 +  dy^2]$, or equivalently the Laplacian has the
normal form $D_s^2 + B(s, y)(y^2 D_s^2 + D_y^2) + LOT$ (lower order terms). We will 
modify Lazutkin's construction
so that in the hyperbolic case the Laplacian has the normal form $D_s^2 + B(s, y)(-y^2 D_s^2 + D_y^2) + LOT$. 
 We do not use the metric
normal form in this paper beyond the quadratic term and we do not really need to define the straightening map
by Lazutkin's method.  We do so anyway because there is no advantage to constructing
another map and because we believe the details of Lazutkin's construction could be useful
in the  general inverse problem.

The straightening map carries $\Delta$ to a variable coefficient Laplacian $\bDelta$
in a neighborhood of $[0, L] \times \{0\}$ in $\R^2.$ A given wave invariant  depends only
on a certain germ of the Laplacian at the orbit,  and consequently only on a certain germ of the
straightening.  So with no loss of generality we may assume that $\bDelta$ is a polynomial
differential operator in the $\y$ variable with analytic coefficients in $\s.$  We define
the Dirichlet Laplacian $\bDelta_{\Omega_o}$ to be this operator with Dirichlet boundary
conditions on $\partial \Omega_o.$  

\subsubsection{Conjugation to normal form}

The next step is roughly to conjugate
$\bDelta_{\Omega_o}$  to a microlocal normal form near $\gamma_o: = [0, L] \times \{0\}$ by an FIO (Fourier integral operator)  on the model space.  Intuitively, we would like to construct a microlocally invertible  FIO $W$ which
'preserves  Dirichlet boundary conditions' and such that $W^{-1} \bDelta W \sim F(|D|, \hat{I})^2$ modulo an
`acceptable remainder' near $\gamma_o$. 

Let us be more precise.  First, by  `preserving  Dirichlet boundary conditions'
we mean that $W$ should carry   the
domain of the Dirichlet Laplacian $\bDelta$ as an unbounded operator  on $\Omega_o$ to the domain of $F (|D|, \hat{I})^2$ as an unbounded operator on $\Omega_o$.  Most significantly, 
$Wu = 0$ on $\partial \Omega_o$ if  $u \in H^1_0(\Omega_o)$. 

Second, let us be more precise about the sense in which we are conjugating $\bDelta$ to $F(|D|, \hat{I})$. It  is technically complicated to conjugate a boundary value problem by Fourier integral operator methods,
so we do something simpler which is sufficient for the proof of our theorem.  Namely we observe that both
 the Dirichlet wave operator $\cos t \sqrt{\Delta_{\Omega}})$ and  the normal form wave operator $\cos t F(|D|, \hat{I})$ 
are restrictions to their domains of well-defined Fourier integral operators in  microlocal neighborhoods in the open space of the bouncing ball orbits . We can use the intertwining operator $W$ on the open space to `pull back'  $\cos t F(|D|, \hat{I})$ (or more precisely its `odd part') to a  parametrix for $\cos t \sqrt{\Delta_{\Omega}})$ in the interior of $\Omega_o$.  It is only in this weak sense that we conjugate the wave group to normal form. 
Since the wave trace $Tr \cos t \sqrt{\Delta}$ 
at $\gamma$ involves the wave kernel only in the interior, we can compute it in terms of this parametrix and hence
in terms of the normal form. 

Our goal then is  the construction of an intertwining operator to normal form in this weak sense  and the characterization of the error term in the conjugation.  There are two main issues we would
like to emphasize in this introduction.  The first has to do with solvability of the conjugation equations. 
Those familiar with the conjugation to normal form in the boundaryless case will recall that (just as in the classical
conjugation to Birkhoff normal form) it is based on solving a  sequence of homological equations $[Q, |D| + \frac{\alpha}{L} \hat{I} ] = \mbox{KNOWN}$ for the infinitesimal intertwining operator $Q$.  This is a first order equation for the symbol of $Q$ and it may seem mysterious that one can solve these equations with two boundary conditions (one at each boundary component).  The second 
point to explain  is the relevant notion  of `acceptable remainder'.  

\subsubsection{Acceptable remainders}

To clarify the second of these points, we recall that the link between the spectrum and normal form is through the coefficients in
the singularity expansion 
\begin{equation} Tr  U(t) = c_{\gamma} (t - L_{\gamma} + i0)^{-1} + a_{\gamma 0} \log (t - L_{\gamma} + i 0) +
\sum_{k = 1}^{\infty}a_{\gamma k}(t - L_{\gamma} + i 0)^k \log (t - L_{\gamma} + i 0) \end{equation}
of the trace of the wave group $U(t) = \exp(it \sqrt{\Delta})$ at $t = L_{\gamma}$. When $L_{\gamma}$
is the length of a bouncing ball orbit (or in general a periodic reflecting ray), 
 the wave trace expansion is very similar to the
boundaryless case in that  the singularity is Lagrangean and the coefficients
may be calculated by the stationary phase method (see Corollary (2.5)).  In another language (cf. Corollary (2.6)), the
wave invariants $a_{\gamma k}$
are non-commutative residues $res (\frac{d}{dt}) U(t) |_{t = L_{\gamma}}$ of the wave group and
its time derivatives.  As already  proved in \cite{G} \cite{Z.1}, it follows that  only a certain amount of data
from the Taylor expansion of the symbol of $\Delta$ along $\gamma$ goes into a given
wave invariant $\gamma.$  The precise statement is that  $a_{\gamma k}$ depends only on the class of $\sigma_{\Delta}$ modulo  the symbol class $S^{2, 2(k + 2)}(V, \R \gamma)$ of Boutet de Monvel \cite{BM}.  Here, $V$ is a conic neighborhood of the symplectic cone $\R \gamma.$ The
bigrading of symbols is in terms of  symbol order and order of vanishing along $\gamma$
(cf. \S 2).  Terms of low symbolic order or of high vanishing order along $\gamma$ do 
not contribute to $a_{\gamma k}$ (cf. Proposition (\ref{DROP})).
Thus we need to construct an FIO $W$ which preserves Dirichlet boundary conditions and which  conjugates $\bDelta_{\Omega_o} \to F(|D_s|, \hat{I})^2$
modulo a remainder in $S^{2, 2(k + 2)}$. 

\subsubsection{Reduction to a semiclassical problem}

 To do this, we convert the problem to a semiclassical
conjugation problem as in \cite{Z.1}. With the proper semiclassical scaling of symbols, elements of $S^{2, K}$ are detected by their coefficient in the semiclassical parameter $N^{-1}.$
The goal then is to construct a semiclassical intertwining operator $W_N$ preserving
Dirichelt boundary conditions and conjugating the scaled version of $\bDelta$ to normal
form modulo a sufficiently high power of $N^{-1}.$

  As in \cite{G}
\cite{Z.1} and elsewhere, we construct $W_N$ as a product $W_N = \Pi_{j = 1}^{\infty} e^{ N^{-j} (P + iQ)_{j/2}}$
as a product of elliptic semiclassical pseudodifferential operators.  The exponents will
just be Weyl pseudodifferential operators $(P + i Q)_{j/2}(s, y, D_y)$ on the transverse space $\R$; here, $P, Q$
are assumed to have real-valued symbols.  The real part $P$ is of two lower orders in $(y, D_y)$ than is
$Q$. As mentioned above, the condition that $W_N$ intertwines to normal form 
 translates into homological equations for the exponents.  
The reason why we can solve these equations while preserving Dirichlet boundary conditions is that  the boundary condition only affects the `odd terms' in $P_{j/2}$ and the  `even terms' in $Q_{j/2}$
with respect to the involution $(y, \eta) \to (y, - \eta).$
The equations for the even/odd parts of $P_{j/2}, Q_{j/2}$ are coupled   except for `diagonal' terms which are
powers of $\hat{I}$.  Hence one can eliminate the `non-diagonal' even parts to get second order ordinary
differential equations for the odd parts.  The boundary problem for these can be solved in the hyperbolic case, and
in the elliptic case they can be solved
as long as $\alpha$ is independent of $\pi$ over $\Q.$  The remaining  powers of $\hat{I}$
constitute the terms in  the normal form.

\subsubsection{Spectral invariance of the normal form }

The intertwining operator $W$ conjugating  $\Delta_{\Omega}$ to $F(|D|, \hat{I})^2$ modulo $S^{2, 2(k + 2)}$ can be
used to construct a microlocal parametrix for the Dirichlet wave group modulo similar acceptable errors. The
definition of the parametrix is ${\cal E}(t) = W F_o(t)\psi_V 1_{\Omega_o}  W_{-1} $, where $W_{-1}$ is a microlocal inverse
to $W$, where $ F_o(t)$ is the odd part of the normal form wave group $\cos t F(|D|, \hat{I})$ (the odd part  satisfying Dirichlet boundary conditions at $s = 0, L$), where $1_{\Omega_o}$ is the characteristic function of $\Omega_o$ and where
$\psi_V$ is a  microlocal cutoff to a suitably small conic neighborhood of $\gamma_o$.  
The wave invariants of the Dirichlet Laplacian $\Delta_{\Omega}$  can then be calculated in terms of the normal form coefficients,  and conversely ( as in \cite{G} (see also \cite{Z.1} \cite{Z.2}) the normal form coefficients can be
determined from the wave invariants.  Hence,  the coefficients of the normal form $F(|D|, \hat{I})$ are spectral invariants of $\Delta_{\Omega}$; a fortiori,  the  classical normal
form of the Poincare map ${\cal P}_{\gamma}$is a spectral invariant.

\subsubsection{Conclusion of the proof}

 Theorem (\ref{TH}) then follows from a theorem  of Colin
de Verdiere \cite{CV} that a bi-axisymmetric analytic plane domain is determined by the
Birkhoff normal form of ${\cal P}_{\gamma}$ of a non-degenerate elliptic  axial orbit $\gamma$.  His proof 
works as well when $\gamma$ is hyperbolic and therefore we can conclude the proof in either case. 

\subsection{Future problems}

It should be remarked that $F(|D|, \hat{I})$ is  explicitly constructed by the algorithm
of this paper.  In conjunction with Lazutkin's construction of a metric normal form, one
can get explicit albeit complicated formulae for the normal form coefficients as polynomials
in the Taylor coefficients of the boundary defining functions at the points $A, B.$  In
the future we plan to take up the obvious question of whether one can determine the
Taylor coefficients from the normal form coefficients.  Colin de Verdiere's theorem shows that the
principal symbol of the normal form alone is enough to determine these coefficients when
the Taylor coefficients at $A$ and $B$ are the same and when the odd coefficients at $A$
vanish.  In less symmetric cases one will have to go into lower order terms in the normal 
form. In the simpler but somewhat analogous case of
surfaces of revolution, it was necessary to go two steps below the principal symbol level
to determine all of  Taylor coefficients from the normal form.  Rotational symmetry is
analogous to left-right symmetry in a plane doman, so we suspect one can solve the
inverse spectral problem at least  for left-right symmetric analytic domains using just the
wave invariants at one orbit.

\subsection{Acknowledgements} The research in this paper was  partially done during visits to the Institut Galilee of Paris Nord in May-June 1998 and to the E.Schrodinger Institute in July 1998.  We also thank M. Zworski for helpful advice on
\S 6 and Y. Colin de Verdiere for asking whether Theorem (\ref{TH}) is also valid in the hyperbolic case.

\section{Preliminaries}

In this section we establish some standard notation and terminology concerning the
wave equation with mixed boundary conditions on a bounded smooth domain and its
associated billiard flow.   For further background we refer to \cite{GM} \cite{PS}.

\subsection{Billiards on plane domains}

First we recall the definition of a non-degenerate bouncing ball orbit for billiards
on a smooth domain in $\R^n$

Let  $\Omega \subset \R^n$ be a  smooth domain, 
 and let $\partial \Omega$ denote its boundary.  Let $g$ be a smooth metric defined
in a neighborhood of $\Omega$. 
 The {\em billiard flow} $\Phi^t$ of
$(\Omega, g)$ is the  geodesic flow  of  $g$ in the
interior $T^*(int \Omega)$,
with the usual law of reflection at the boundary  (i.e.
the tangential component of the velocity remains the same but the normal component changes
its sign). 
Here, $int \Omega$ denotes the interior of $\Omega.$

From the symplectic point of view, there are two important hypersurfaces in
$T^*(\Omega)$: the manifold $T^*_{\partial \Omega} \Omega$ of (co)vectors with
footpoints on $\partial \Omega$ and the unit cotangent bundle $S^*\Omega$. We
understand by $T^*\Omega$ the restriction $T^* \R^2|_{\Omega}$. The intersection
$S^*_{\partial \Omega}\Omega = T^*_{\partial \Omega} \Omega \cap  S^* \Omega$ is transversal. 
The characteristic foliation of $S^*\Omega$ is spanned by the generator of the geodesic
flow.  The  projection $\pi : T^*_{\partial \Omega} \Omega  \rightarrow T^*\Omega$ which
projects a (co)vector at $x \in \partial \Omega$ to the (co)tangent hyperplane defines
a real line whose fiber at $x$ is the normal bundle $N_x (\partial \Omega)$.  We denote
by $\nu_x$ the inward unit normal so that $N_x (\partial \Omega) = \R \nu_x.$ This
line bundle coincides with the
 characteristic foliation of $T^*_{\partial \Omega} \Omega$. The image of $S^*_{\partial
\Omega} \Omega$ under $\pi$ is the unit disc bundle $D^*(\partial \Omega).$ We will
identify it with the space $S^{in}_{\partial \Omega}\Omega$ of inward unit tangent
(co)vectors with footpoints on $\partial \Omega.$ 

 One then defines the {\em billiard ball map}
$\beta: D^*(\partial \Omega) \to D^*(\partial \Omega)$ as follows: lift a tangent (co)vector
$v$ to $\partial \Omega$ of length $<1$ to $S^{in}_{\partial \Omega}\Omega$ and move
along the corresponding geodesic until the ball hits the boundary.  Then project its tangent
vector to $D^* \partial \Omega.$

By a {\em reflecting ray} one means a broken geodesic or billiard trajectory of the
billiard flow whose intersections with the boundary are all transversal. We denote
the successive points of contact with the boundary by $q_0, q_1, q_2, \dots.$  Of special
importance here are the {\em periodic reflecting rays} where $q_n = q_0$ for some $n > 1.$
We will denote such a periodic trajectory by $\gamma : [O, T] \rightarrow S^*\Omega$.  By a {\em bouncing ball} orbit $\gamma$ one means a periodic reflecting ray where $q_0 = q_2.$
The projection to $\Omega$ consists of a segment $q_0q_1$ which is orthogonal to the
boundary at both endpoints, i.e. $q_0q_1$ is an extremal diameter.  The period is of course twice the length of the segment,
which we denote by $L$.  We write $q_0 = A, q_1 = B.$

The Poincare map ${\cal P}_{\gamma}$ of a periodic reflecting ray $\gamma$ is the first return map
to a symplectic transversal.  A natural symplectic transversal is given by a neighborhood of
 $\gamma(0)$ in 
$S_{\partial \Omega}^{in}\Omega$. Thus for
$v \in S_{\partial \Omega}^{in}\Omega$,  ${\cal P}_{\gamma}(v)$ is obtained by following
the broken geodesic thru $v$ until it reflects from $\partial \Omega$ the $n$th time in an 
inward vector. The linear Poincare map $P_{\gamma}$ is defined to be $d {\cal P}_{\gamma}(\gamma(0)).$

Equivalently, one can idenfity $S_{\partial \Omega}^{in}\Omega$  with $D^*(\Omega)$ as
above and define ${\cal P}_{\gamma}$ as a map on a neighborhood $S_{\gamma}$ of the zero vector at an endpoint of $\gamma$ in $D^*(\Omega)$.  Then for a periodic reflecting ray $\gamma$ with $n$
reflections, ${\cal P}_{\gamma}$ may be identified with $\beta^n |_{S_{\gamma}}.$

\begin{defn} A bouncing ball orbit $\gamma$ is said to be non-degenerate if the eigenvalues of $P_{\gamma}$ are not
roots of unity.  There are two cases: $\gamma$ is 

\noindent(i)  {\em non-degenerate elliptic} if the eigenvalues
of $P_{\gamma}$ are of the form $\{e^{\pm i \alpha}\}$ with $\alpha/\pi \notin \Q.$;

\noindent(ii) {\em non-degenerate hyperbolic}  if the eigenvalues
of $P_{\gamma}$ are of the form $\{e^{\pm \lambda}\}$ for some $\lambda \in \R^+$. 
\end{defn}

For background on bouncing-ball orbits we refer to \cite{B.B} (see \S 5.2). If we denote by $R_A$, resp. $R_B$
the radii of curvature of the boundary at the endpoints $A$, resp. $B$ of an extremal diameter, then one finds
that ellipticity is equivalent to the condition that  $L < R_A + R_B, L > R_A, L>R_B$ or else that
$L < R_A + R_B, L < R_A, L < R_B$.  If $\bar{AB}$ is a local minimum diameter than $L < R_A + R_B$ while if it
is a local maximum diameter then $L < R_A + R_B$.  So a (non-degenerate) local maximum diameter must be hyperbolic;
 a local minimum diameter is elliptic if it satisfies the additional inequalities above.  Under our symmetry assumption
$R_A = R_B$,  so ellipticity is equivalent to the statement that the center of curvature at $A$ lies below the horizontal
axis.

\subsection{Wave trace on a manifold with boundary}

Let $\Delta$ denote the Laplacian of a metric $g$ defined in an open neighborhood of a bounded smooth
domain $\Omega \subset \R^2$.   The Dirichlet Laplacian $\Delta_{\Omega}$ of $\Omega$ is then defined to
be the self-adjoint operator on $L^2(\Omega)$ with domain 
\begin{equation} {\cal D}om (\Delta_{\Omega}) = \{u \in H_0^1(\Omega) : \Delta u \in L^2(\Omega)\} \end{equation}
where $H_0^1(\Omega)$ denotes 
the closure of $C_0^{\infty}(\Omega)$ under the norm
$$||u||_1^2 = \sum_{j \leq 2} ||\frac{\partial u}{\partial x_{j}}||^2_{L^2}. $$
and where 
where $\Delta u$ is taken in the sense of distributions.  We recall that $H_0^1(\Omega) = \{u \in H^1(\R^2): \mbox{supp}\;u
\subset \Omega\}.$

Now let $E(t,x,y)$ be the fundamental solution of the mixed wave equation with Dirichlet
boundary conditions:
\begin{equation}\label{DWK} \begin{array}{ll} \frac{\partial^2 E}{\partial t^2} =
\Delta E & \mbox{on} \; \R \times \Omega \times \Omega \\ & 
\\ E(0, x,y) = \delta(x - y) & \frac{\partial E}{\partial t}(0,x,y) = 0 \\
& \\
E(t,x,y) = 0 & (t, x, y) \in \R \times \partial \Omega \times \Omega. \end{array} \end{equation}

As discussed in (\cite{GM}, \S 5]\cite{PS}), the fundamental solution is not globally a Lagrangean
distribution.  However, for any $T > 0$ there exists a  conic neighborhood $\Gamma_T$ of the bouncing ball orbit $\gamma$ so that the microlocalization of  $E(t, x, y)$ to
$\Gamma_T$  is Lagrangean for $|t| \leq T.$  More
precisely, 
\begin{theo}\label{POISSON} (\cite{GM}, Theorem 4.1 and Proposition; or \cite{PS}, \S 6) Let $(x, \xi) \in T^*(int \Omega)$ or let $x \in \partial \Omega$ and suppose that $\xi$ is a non-glancing (co-) direction at
$x$.  Then there exists a conic neighbhorhood ${\cal O}$ of $(x, \xi)$ in $T^*\R^n - 0$
and a Fourier integral distribution $\tilde{V}(T, x, y)$ essentially supported in
${\cal O}$ such that for any pseudodifferential operator $\chi_{{\cal O}}(x, D)$
essentially supported in ${\cal O}$, $\chi_{{\cal O}}(x, D) \tilde{V}(T, x, y) -
\chi_{{\cal O}}(x, D)E(T, x, y) \in C^{\infty}(\R \times \Omega \times \Omega).$
\end{theo}

By a partition of unity, this leads to a Fourier integral formula for the wave trace near a periodic
reflecting ray:

\begin{theo} \cite{A.M} Suppose that $T$ is an isolated point of the length spectrum of
$\Omega$ with the following properties: 

\noindent(i) $T \not= |\partial \Omega|$;\\
(ii) All of the closed billiard trajectories of length $T$ are non-degenerate
reflecting rays.

Then modulo smooth functions in $t$ near $t = T$ we have:
$$ \begin{array}{l} \sum_{\lambda_j \in Sp({\sqrt{\Delta}})} cos \lambda_j t =
\int_{\Omega} E(t,x,x) dx  \\ \\
= \sum_j  \int_{\Omega} \tilde{V}_{\pm}^{(j)}
(t,x,x) dx \end{array}$$where $\tilde{V}_{\pm}^{(j)}(t,x,y)$ are Fourier integral distributions 
associated to the broken geodesic flow near the closed trajectories of length $T$. \end{theo} 

The wave invariants can therefore be obtained  by applying the method
of stationary phase to an oscillatory integral.  An apparent obstruction is
that the domain of integration is a manifold with boundary. However, Guillemin-
Melrose prove:

\begin{lem}\label{MSP} Let $\Omega \subset \R^n$ be a smooth domain with boundary and
let $\phi \in C^{\infty}(\R^n)$ have clean critical point sets.  For each
critical value $\lambda$ let $C_{\lambda}$ denote the critical set in $\R^n.$
Suppose that for each $\lambda$, $C_{\lambda}$ intersects $\partial \Omega$
transversally. Let $U_{\lambda}$ be a neighborhood of $C_{\lambda}$ with
the properties:

\noindent(i) $U_{\lambda}$ contains no critical points of $\phi$ except
$C_{\lambda}$;\\
(ii) $U_{\lambda} \cap \partial \Omega$ contains no critical points of
$\phi|_{\partial \Omega}$ except $C_{\lambda} \cap \partial \Omega.$ Then
the usual stationary phase expansion is valid:
$$\int_{\Omega} a(x) e^{i \tau \phi(x)} dx \sim e^{i \tau \lambda} \tau^{-k/2}
\sum_{i = 0}^{\infty} \alpha_i \tau^{-i},\;\;\;\;\;(\tau \rightarrow \infty)$$
where  $k$ is the codimension of $C_{\lambda}$ in $\R^n$ and where $\alpha_i$
are integrals of $a$ and its derivatives over $C_{\lambda} \cap \Omega.$
\end{lem}

\subsection{Wave trace invariants as non-commutative residues}

We may summarize the relevant result on the Poisson formula as follows:

\begin{cor} Let $\gamma$ be a non-degenerate billiard trajectory whose length $L_{\gamma}$
is isolated and of multiplicity one in $Lsp (\Omega)$.  Let $\Gamma_L$ be a sufficiently small conic neigbhorhood
of $\R^+ \gamma$ and let $\psi$ be a microlocal cutoff to $\Gamma_L$.  Then for
$t$ near $L_{\gamma}$, the trace of the wave group has the singularity expansion 
$$ Tr \psi E(t) = c_{\gamma} (t - L_{\gamma} + i0)^{-1} + a_{\gamma 0} \log (t - L_{\gamma} + i 0) +
\sum_{k = 1}^{\infty}a_{\gamma k}(t - L_{\gamma} + i 0)^k \log (t - L_{\gamma} + i 0)$$
where the coefficients $a_{\gamma k}$ are calculated by the stationary phase method
from a Lagrangean parametrix. \end{cor}

This corollary allows us to identify the wave invariants as non-commutative residues as in \cite{G}\cite{Z.1}.  Recall that if $A$ is a Fourier integral operator, and if $P$
is any positive elliptic first order pseuodifferential operator, then the zeta function
$\zeta(z, A, P): = Tr A P^{-z}$ has  meromorphic extension to $\C$ with at most simple
poles.  The residue at $z=0$ is referred to as the non-commutative residue $res (A)$
of $A$.  It is independent of $A$ and is a tracial invariant, i.e. $res(W A W^{-1})=
res A.$  In the boundaryless case, one has $a_{\gamma k} = res (\frac{d}{dt})^k E(t) |_{t = L_{\gamma}}$.  The only ingredients in the proof are the Lagrangean property of $E$
and a canonical transform between $Tr A e^{i t P}$ and $Tr A P^{-s}.$ Hence the same
result remains valid in the case of periodic reflecting rays of the boundary case:

\begin{cor} \label{RES} If $\gamma$ is a periodic reflecting ray, then
$a_{\gamma k} = res (\frac{d}{dt})^k E(t) |_{t = L_{\gamma}}$. \end{cor}

Now we recall some elementary results from \cite{G} and \cite{Z.1} on the data of the domain or metric
which go into a given wave invariant $a_{\gamma k}$ in the boundaryless case. Analogous results hold in
the boundary case, but we postpone stating them.

In the following, let $P$ be any first order pseudodifferential operator of real principal type on a boundaryless manifold
$U$
and assume for simplicity that all closed orbits of its bicharacteristic flow are non-degenerate.  
Let  $a_{\gamma k}(P)$ denote the kth wave invariant at a closed orbit $\gamma$ for $e^{it P}$. 
Since it is calculated by a stationary phase expansion at $\gamma$, it is obvious that  $a_{\gamma k}(P)$ depends only on certain part of the jet of the complete symbol of $P$
around $\gamma$.  

 To state the precise result, let us Taylor  expand each term in   the complete (Weyl) symbol $p(s, \sigma, y, \eta) \sim p_1 + p_o + \dots $ of $P$ at $\R^+\gamma$: 
$$p_j(s, \sigma, y, \eta) = \sigma^j
 p_j(s, 1, y, \frac{\eta}{\sigma})= \sigma^j (p_j^{[o]} + p_j^{[1]} + \dots)$$
with 
$p_j^{[m]}(s, 1, y, \frac{\eta}{\sigma})$ the part which is homogeneous of
 degree $m$ in $(y, \frac{\eta}{\sigma}).$  Set $P_j:=Op(p_j), P_j^{[m]} := Op(p_j^{[m]})$
and $P^{\leq N}_j = \sum_{m\leq N} P_j^{[m]}$. Then we have:
\begin{prop} \label{DROP} (\cite{Z.1}, Proposition 4.2)
$a_{\gamma k}(P) = a_{\gamma k}(P_1^{\leq 2(k+2)}  + P_o^{\leq 2(k+1)}
 + \dots P_{-k -1}^{o}).$ \end{prop}

Thus, it is sufficient to define the normal form and the intertwining operator  in a microlocal (conic) neighborhood $V$
of $\R^+ \gamma$.  In $(s, \sigma, y, \eta)$ coordinates, we may define the cone by:
\begin{equation} V = \{(s, \sigma, y, \eta) \in T^*U_{\epsilon} : |y| \leq \epsilon, |\eta| \leq \epsilon |\sigma|\}.
\end{equation} 
 Let $\psi_{V}$ be a microlocal cutoff to $V$, with symbol identically equal
to one in a slightly smaller open cone around $\R \gamma$ and put $\Delta_{\epsilon}:= \psi_{\epsilon} \Delta_{\Omega}
\psi_{\epsilon}$.  $\Delta_{V}$ and  $\Delta$ of course have the same microlocal normal form around
$\gamma$ so for      notational simplicity we  often drop the subscript and 
leave it the reader to recall that the operator is cutoff.  We will re-instate the cutoff at the crucial point of
calculating the residue.

Since terms which vanish to too high order at $\gamma$ or which have too low a pseudodifferential order do not
contribute to $a_{\gamma k}$ we introduce   a bi-grading on symbols in terms of order as
a symbol and order of vanishing along $\gamma$.   
 Let $\R^+ \gamma \subset T^* U$ be the cone thru an embedded curve $\gamma \in T^* U$ and $O_j S^m(U)$  denote the class of symbols of order $m$  over $U$ which vanish to order $j$ along $\gamma$.  Following  (\cite{BM}, see also
\cite{BMGH})
 we denote by $S^{m, k}(V, \R^+ \gamma)$
the  the class of symbols microsupported in $V$ 
 which admit asymptotic expansions
\begin{equation}\label{DGS} a \sim \sum_{j = 0}^{\infty} a_{m - j}, \;\;\;\;
a_{m - j}(x, r \xi) = r^{m - j} a(x, \xi),\;\;\;\;\; a_{m - j} \in O_{k - 2j}S^{m - j}.
\end{equation}
Here,  $k \in {\bf N}$ and is no condition if $2j \geq k.$  We denote by $Op^w S^{m,k}$ the Weyl pseudodifferential
operators with complete symbols in $S^{m,k}.$ 

In local coordinates $(s, y)$ with $\gamma = \{y = 0\}$ and with dual symplectic coordinates $(\sigma, \eta)$,
  $a \in  S^{m, k}(V, \R^+ \gamma)$ if
 \begin{equation}\label{HOMO}\begin{array}{l} a(s, \sigma, y, \eta)
 \sim \sum_{j = 0}^{\infty} a_{m - j, k}(s, \sigma, y, \eta) 
\sim \sum_{j = 0}^{\infty} \sigma^{m - j} a_{m - j, k}(s, 1, y, \eta/\sigma)\\ \\
a_{m - j, k}(s, 1, y, \eta/\sigma) \in O_{k - 2j}S^0.
\end{array} \end{equation} 
The coefficient $a_{m - j}(s, 1, y, \eta/\sigma)$ is homogeneous of degree $0$ and
the assumption that $a$ is microsupported in $V$ becomes that  $a_{m - j}$
is supported in the transverse ball $B_{\epsilon} = \{|(y', \eta')| < \epsilon.\}$

Any symbol $a \in S^m(V)$ may be expanded as a sum of symbols in $S^{m ,k}$.  Indeed, let $a_{m - j; k}$ be the term  of degree $k - 2j$ in its Taylor expansion for $k \geq 2 j$
\begin{equation}\label{TAY} a_{m - j; k } =  \sum_{|\alpha| + |\beta| = k - 2j} a_{m - j; k - 2j,
\alpha, \beta}(s, ) y^{\alpha} (\eta / \sigma)^{\beta}. \end{equation}
  Then
the asymptotic sum $a_{(m, k)} \sim \sum_{j = 0}^{[k/2]} a_{m - j; k  }$ belongs to
$S^{m, k}$ (sharp) and $a = \sum_{k = 0}^{\infty} a_{(m, k)}.$  In particular, we may
expand $\Delta$ in this form, and $a_{\gamma k}$ depends only on the class of $\Delta$
modulo $S^{2, k}.$

We summarize the discussion by restating Proposition (\ref{DROP}) in terms of these symbol classes:
 \begin{cor} Suppose that $A, B \in Op^w S^1(V)$ and that $A \equiv B \;\mbox{mod}\; Op S^{1, k}(V)$.  Then
$a_{\gamma k}(A) = a_{\gamma k}(B).$  
 \end{cor}

\section{Straightening the domain}

As discussed in the introduction, the normal form lives on the model domain
$\Omega_o = [0, L] \times \R$, or more precisely in a microlocal neighborhood of $T^*[0, L] $
in its cotangent bundle. In this section we introduce the analytic objects on the model
domain and explain how to transfer $\Delta$, in a 
 neighborhood in $\R^2$ of a non-degenerate (elliptic or hyperbolic)
bouncing ball orbit $\gamma$, to a variable coefficient Laplacian $\bDelta$ in  a neighborhood of $[0, L] \times \{0\}$ in
the model.   We emphasize that all maps and operators that we discuss in this section extend
to  open domains containing the various manifolds with boundary.

\subsection{The model domain $\Omega_o$}

The  configuration space of the model is the infinite strip $ \Omega_o$. We denote the coordinate on $[0,L]$ by  $s$  and that on $\R$ by  $y$ , with dual cotangent
coordinates $\sigma, \eta$ on $T^* [0,L] \times T^* \R$.  We also denote by $I_e = \half(\eta^2 + y^2)$ the elliptic action variable and by $I_h = \half(\eta^2 - y^2)$ the hyperbolic action variable 
on $T^*\R$. To simplify notation we often just write $I$ for the relevant action variable.
In the Poisson algebra of $T^*(\R^2)$ we consider the maximal abelian
subalgebra   ${\cal A}_{cl} = \langle |\sigma|, I \rangle.$  The model (classical)
Hamiltonians are those of the form $H_{\alpha} = |\sigma| + \frac{\alpha}{L} I$ which generate
 linear Hamiltonian flows. As in the case of straightline motion, $H_{\alpha}$ 
generates a broken
Hamiltonian flow on $T^* [0,L] \times T^* \R$ when equipped with the boundary condition that the trajectory is reflected by
$\tau$ when it hits the boundary.

We will view $\Omega_o$ as a submanifold with boundary of the  
`open space' $S^1_{2L} \times \R$
with $S^1_{2L} = \R / 2 L \Z \cong [-L, L]$. Many of our operators will be supported  in a  neighborhood 
$U_{\epsilon} = (- \epsilon, L + \epsilon)
\times \R)$ of $\Omega_o$ in $S^1_{2L} \times \R$.  We denote by $\tau_{\epsilon}(s)$
a smooth cutoff to $U_{\epsilon}$ with $\tau_{\epsilon} \equiv 1$ in a smaller neighborhood
of $\Omega_o.$

Since $\Omega_o$,  $U_{\epsilon}$, and $S^1_{2L} \times \R$ are products,  their algebras of pseudodifferentials
operators are easily described in terms of  pseudodifferential operators along $S^1_{2L}$ and 
transverse pseudodifferential operators on $\R.$ 

In the direction of $S^1_{2L}$, we introduce 
the algebra $\Psi^*(S^1_{2L}) = \langle s, D_s \rangle$ of standard pseudodifferential operators on $S^1_{2L}$.  We distinguish the element $|D_s| = \sqrt{- D_s^* D_s}$ with eigenfunctions
$e^{i \frac{\pi k s}{L}}$  associated eigenvalues $|k|$. 
We also define the subspaces and projections:
\begin{equation}\begin{array}{ll} H^2_{\pm} (S^1_{2L}) = \oplus_{k = 0}^{\infty} \C e^{\pm i \frac{\pi k \s}{L}} & \Pi_{\pm} : L^2
\to H^2_{\pm} \\ & \\
\Pi_{+k}: L^2(S^1_{2L}) \to  \C e^{ i \frac{\pi k \s}{L}} & \Pi_{- k}: L^2(S^1_{2L}) \to  \C e^{-  i \frac{\pi k \s}{L}}
\\ & \\
L^2_{odd} (S^1_{2L}) = \oplus_{k = 0}^{\infty} \C \sin( \frac{\pi k \s}{L}) & \Pi_o : L^2 \to L^2_{odd},\;\;\;
\Pi_o f(s) = \half ( f(s) - r f(s))\\ & \\
 \Pi_{ok}:  L^2 (S^1_{2L}) \to \C \sin( \frac{\pi k \s}{L}) & 
\end{array} \end{equation} 
Here, $r f(s) = f(-s)$ where $-s$ is taken modulo $2L$, i.e.  $r$ is reflection through
the boundary. 
We may tensor the subspaces with $L^2(\R)$ to get corresponding subspaces of and operators on  $L^2(S^1_{2L} \times \R)$ and we use the same notation for these.

We will sometimes identify functions on $\Omega_o$ with odd functions on $S^1_{2L} \times \R.$ 
More precisely, let us put:
\begin{equation}\begin{array}{ll} 
\Pi_k: L^2([0, L] \times \R) \to \C \sin( \frac{\pi k \s}{L})\otimes L^2(\R) &
\\ & \\ 
1_{[0,L]}(\s) : L^2(S^1_{2L}) \to L^2([0, L]) & 1_{[0,L]}(\s) = \mbox{the characteristic function of}\;
[0, L] \\ & \\
A: L^2([0, L]) \to L^2_{odd}(S^1_{2L}) & Af =  \;\mbox{the odd extension of } \; f\\ & \\
A^* : L^2_{odd}(S^1_{2L}) \to L^2([0, L]) & A^* g = 1_{[0, L]} \Pi_o g \\ & \\
|D_{s}| : H^1_0([0, L]) \to L^2([0, L]) & |D_s| \sin( \frac{\pi k \s}{L}) = \frac{\pi k}{L} 
\sin( \frac{\pi k \s}{L})  \end{array}
\end{equation}

We use the  notation $|D_s|$ (or simply $|D|$) simeltaneously for the `Laplacian' on $S^1_{2L}$ and
the Dirichlet Laplacian on $[0, L]$. No confusion should result since they are
defined on different domains, and moreover the definitions are compatible under the
above identification, as the following proposition shows. 
\begin{prop} We have:

\noindent(i) $A : L^2([0, L]) \to L^2_{odd}(S^1_{2L})$ and  $A A^* = \Pi_o, A^* A = Id.$ 

\noindent(ii) $A^* |D_s| A = |D_s|$

\end{prop} 

The proof is obvious so we omit it.

In the transverse direction, we first introduce the isotropic Weyl algebra
${\cal W}^*$.  This is a completion of  the algebra  
${\cal E} := <y, D_y> $ 
of polynomial differential operators on $\R$. We
denote by 
${\cal E}^n$ denote the subspace of polynomial differential operators of degree n in the variables $y, D_y.$
We also denote by ${\cal E}^n_{\epsilon}$ the polynomials all of whose terms have the same parity as $n$. 
In the isotropic Weyl algebra ${\cal W}^*$, the
operators $y, D_y$ are given the order $\frac{1}{2}$, so that
\begin{equation} {\cal E}^n \subset {\cal W}^{n/2}, \;\;\;\;
[{\cal E}^m, {\cal E}^n] \subset {\cal E}^{m + n - 2}. \end{equation}
We note that such operators are not standard homogeneous pseudodifferential operators, but can be rescaled to this form. The rescaled algebra is generated by the first order
homogeneous pseudodifferential operators $y |D_s|^{\half}$ and $ |D_s|^{-\half} D_y.$
In the open space the relevant algebra is the algebra of homogeneous pseudodifferential operators over $U_{\epsilon}$ generated by $D_s, y |D_s|^{\half}, |D_s|^{-\half} D_y.$ 
To be more precise, this construction is only well-defined in a microlocal neighborhood
of $\gamma$ where $|D_s|$ is elliptic.

It should also be recalled that $a^w(y, D_y) \in {\cal E}^n$ with real-valued symbols are essentially self-adjoint
operators on ${\cal S}(\R)$ (the Schwartz space), i.e. have a unique self-adjoint extension to $L^2(\R).$ Hence
their exponentials $e^{i a^w(y, D_y)}$ are unambiguously defined.

 \subsection{Half-density Laplacian}

In order to deal with self-adjoint operators with respect to the Lebesgue density $ds dy$, we 
pass from the scalar Laplacian to the (unitarily equivalent) 1/2-density
Laplacian
\begin{equation} \begin{array}{l} \Delta_{\half} := J^{1/2} \Delta J^{-1/2}\\
\\ = \sum_{i,j = 1}^2 J^{-1/2} D_{x_i} g^{ij}J D_{x_j} J^{-1/2} \\
\\
=  g^{11} D_s^2 + g^{22} D_y^2 + 2 g^{12} D_s D_y + \frac{1}{i} \Gamma^1
D_s + \frac{1}{i} \Gamma^2 D_y + \sigma_0.\end{array} \end{equation} 
Here we write $(s,y) = (x_1, x_2), D_{x_j} = \frac{\partial}{i \partial x_j}.$ 
The functions $\Gamma_j$ are real valued.  Since it is self-adjoint relative to
the Lebesgue density, 
its complete Weyl symbol is real valued.  
 Henceforth we denote the 1/2-density Laplacian simply by $\Delta.$

\subsection{Straightening the domain}

Let us now explain how to transfer the Laplacian to the model domain. 
  As in the boundaryless case,  the model
space is in some sense the normal bundle of the orbit.
This is literally correct in the boundaryless case and  the exponential map
along the normal bundle $N_{\gamma}$ of $\gamma$ can be used to transfer $\Delta$ 
 to the normal bundle.  In the case of
a bouncing ball orbit in the  boundary case, the normal bundle and exponential map are ill-defined at the reflection points  
but  Lazutkin has  constructed a nice replacement for them.  Namely, he constructs a map $\Phi$ which straightens
the domain to a strip near a stable  elliptic bouncing ball orbit and which simeltaneously puts the Laplacian into a preliminary normal form.   We will modify his method to encompass hyperbolic bouncing ball orbits as well. 
 As mentioned above, we do not need the full details of the map or metric
normal form here.  Hence we only sketch the construction of $\Phi$, referring the reader
to \cite{L} for the details. In the following, $\Omega_{\epsilon}$ denotes an $\epsilon$-neighborhood of $\overline{AB}$ in 
$\Omega$.  

\begin{defn} By a transversal power series map from $\Omega_{\epsilon}$
to $\Omega_0$ we mean a formal power series 
$$ \Phi : \Omega_{\epsilon}
\rightarrow U_{\epsilon},\;\;\;\;\; \Phi(s, y) = (\s, \y) $$
of the form
$$ \begin{array}{l} \s = s + \sum_{m = 2}^{\infty} 
\kappa_{m m} (s)   y^m \\ \\ 
\y = \sum_{p=1}^{\infty}  \psi_{p p}(s)   y^p \end{array}$$
with  real valued analytic  coefficients in $s$ extending analytically to a neighborhood
of $[0, L]$ in $\C$  and
satisfying the boundary conditions:
$$ \{\s = 0\} \cup
\{\s = L\} = \Phi (\partial \Omega_{\epsilon}). $$\end{defn}

We use the language of formal power series since we only need to use a polynomial part
of the map to construct the normal form up to a desired accuracy. 
 Indeed, 
a given  wave invariant $a_{\gamma k}$ only involves the $2k + 4$ -jet  of $\Phi.$  The convergence
of the series is irrelevant to our purposes and we will not discuss it. 
 The following
Lemma was in effect proved by Lazutkin in \cite{L} in the elliptic case.

\begin{lem} Suppose that $\overline{AB}$ is a bouncing ball orbit. Then there exists a transversal power series  map $\Phi : (\Omega_{\epsilon}, \partial
\Omega_{\epsilon}) 
\rightarrow (\Omega_o, \partial \Omega_o)$ which straightens the
domain and puts  $\bDelta$ in the form:
$$\begin{array}{ll} \bar{\Delta} \sim D_{\s}^2 + B(\s, \y)(\y^2 D_{\s}^2 + D_{\y}^2) +
\Gamma'_s D_s + \Gamma'_y D_y & \rm{elliptic \; case} \\ & \\
\bar{\Delta} \sim D_{\s}^2 + B(\s, \y)( D_{\y}^2 - \y^2 D_{\s}^2) +
\Gamma'_s D_s + \Gamma'_y D_y & \rm{hyperbolic \; case} \end{array}
$$
in the sense that the left and right sides agree to infinite order at $y = 0.$
Here, $B(s,y)$ is a transversal power series.
\end{lem}

Of course,  we only need this form of $\bDelta$ to order $K$
construct the normal form modulo  $Op S^{2, K}$, and only in the principal terms
do we need to know the exact form. Therefore we  only briefly recall the proof of the
lemma and  only discuss the principal terms in detail.

 Under any map $\Phi$, the usual (scalar) $\Delta$ conjugates to
  \begin{equation}\bDelta_0:= \Phi^{* -1}\Delta \Phi^* =  |\nabla \s|^2
D_{\s}^2 + |\nabla \y|^2 D_{\y}^2 +   2\langle \nabla \s, \nabla \y\rangle
D_{\s} D_{\y} + \frac{1}{i} \Delta \s D_{\s} +\frac{1}{i} \Delta \y D_{\y}.
\end{equation}
In the case of a transversal power series map, the coefficients 
are also transversal power series of the form: 
\begin{equation}\label{PS} \begin{array}{l} |\nabla \s|^2 = 1 + 
\sum_{m=2}^{\infty} a_{m m}(\s) \y^m \\ \\ 
|\nabla \y|^2 = \sum_{m=0}^{\infty} b_{m m}(\s) \y^m\\ \\ 
\Delta \s = \sum_{m=0}^{\infty} c_{m m}(\s) \y^m\\ \\ 
\Delta \y = \sum_{m=0}^{\infty} d_{m m}(\s) \y^m \\ \\
\langle \nabla \s, \nabla \y \rangle = \sum_{m=1}^{\infty}
e_{m m}(\s) \y^m    \end{array} \end{equation} 
The volume density in the new coordinates,  $J(\bar{s}, \bar{y})$ has a similar form.
The 1/2-density Laplacian $ (\Phi^{* -1}\Delta \Phi^*)_{\half}$ in the
transformed coordinates is   given by
  \begin{equation}\bar{\Delta}:= J^{-\half} (\Phi^{* -1}\Delta \Phi_0^*)
J^{\half} =  |\nabla \s|^2 D_{\s}^2 + |\nabla \y|^2 D_{\y}^2 +   2\langle \nabla
\s, \nabla \y\rangle D_{\s} D_{\y} +  \frac{1}{i} \Gamma_{\s}  D_{\s} + \frac{1}{i} \Gamma_{\y}  D_{\y}
+ K \end{equation} where
\begin{equation}
 \begin{array}{l}
 \Gamma_{\s} = - |\nabla \s|^2
\frac{\partial}{\partial \s} \log J - \langle \nabla \s, \nabla \y \rangle
\frac{\partial}{\partial \y} \log J +     \Delta \s  \\ \\ \Gamma_{\y} = - |\nabla
\y|^2 \frac{\partial}{\partial \y} \log J - \langle \nabla \s, \nabla \y \rangle
\frac{\partial}{\partial \s} \log J + \Delta \y  \\ \\ K =
J^{-\half}(\Phi^{* -1}\Delta \Phi^*) (J^{\half}).  \end{array} \end{equation}
Combining with the power series expressions in (\ref{PS}) one has an expression
for $\bDelta$ modulo terms vanishing to order $K$ at $y = 0.$

The principal term (modulo  $S^{2, 3}$) is given by $D_{\s}^2 + (a_{22}(s)
\y^2 D_{\s}^2 + e_{11}(s) \y D_{\y} + b_{00}(s) D_{\y}^2.$  
We now choose the coefficients in $\Phi$ to put it in the preliminary normal form of
\begin{equation}\begin{array}{ll}   D_{\s}^2  +      \half \dot{b}_{00}(s)[\y^2 D_{\s}^2 +  D_{\y}^2] & \rm{elliptic \;
case} \\ & \\D_{\s}^2  +      \half \dot{b}_{00}(s)[  D_{\y}^2 - \y^2 D_{\s}^2 ] & \rm{hyperbolic \;
case}\end{array}
.
\end{equation}

This requires the coefficients to solve the equations:
\begin{equation} \begin{array}{lll} \bullet & e_{11} = 0 &  \\ & & \\
\bullet & a_{22} =  b_{00} & \rm{elliptic \; case} \\ & & \\
\bullet & a_{22} = -   b_{00} & \rm{hyperbolic\; case}
.\end{array} \end{equation}

To solve these equations we use the following expressions from (\cite{L}, (4.5) -- (4.9)):
 \begin{equation}\begin{array}{l}
a_{22}(s) = (2 \kappa_{22}' + 4 \kappa_{22}^2) \psi_{11}^{-2} \\ \\
b_{00}(s) = \psi_{11}^2(s) \\ \\
c_{00}(s) = 2 \kappa_{22}(s) \end{array} \end{equation}
Since they are obtained purely algebraically, the same equations hold in both the elliptic and hyperbolic cases.
They are easiest to solve if we make the substitution $\theta_{11} = \psi_{11}^{-1}.$
One then has:
\begin{equation}\begin{array}{l} e_{11} = 0  \Rightarrow \kappa_{22}(s) = - \half \frac{\psi_{11}'(s)}{\psi_{11}(s)} 
\\ \\
a_{22} = b_{00} \Rightarrow 
\theta_{11}''  = \frac{1}{\theta^3}\;\;\; \rm{elliptic \; case} \\ \\ 
a_{22} = - b_{00} \Rightarrow 
\theta_{11}''  = -  \frac{1}{\theta^3}\;\;\; \rm{hyperbolic \; case}\end{array} \end{equation} 
The solutions have the form:
\begin{equation}\begin{array}{ll}  \theta_{11} = \sqrt{\ell +
\frac{(\s - s_0)^2}{\ell}} & \rm{elliptic \; case} \\ & \\
 \theta_{11} = \sqrt{- \ell +
\frac{(\s - s_0)^2}{\ell}} & \rm{hyperbolic \; case} \end{array}
\end{equation}
for some constants $\ell, s_0$. 
To verify this it is easiest to substitute $\xi = \theta_{11}^2$ . The equation for $\xi$ is then
$\half \xi \xi'' = \frac{1}{4} (\xi')^2 \pm 1$ (+ = elliptic, - = hyperbolic).  
Solving for $b_{00}$ gives 
$$  b_{00}(\s) = 
 \frac{1}{\ell + \frac{(\s - s_0)^2}{\ell}}\;\; (\rm{elliptic}), \;\;\;\;\;\; b_{00}(\s) = 
 \frac{1}{- \ell + \frac{(\s - s_0)^2}{\ell}}\;\;\; (\rm{hyperbolic})$$
hence
\begin{equation} \int_0^{\s} b_{00}(\s) d\s = \frac{1}{\ell} \tan^{-1} (\frac{\s -
s_o}{\ell}) \;\; (\rm{elliptic}), \;\;\;\;\;\;\int_0^{\s} b_{00}(\s) d\s = \frac{- 1}{\ell} \tanh^{-1} (\frac{\s -
s_o}{\ell})\;\;\; (\rm{hyperbolic}) . \end{equation}
  In the elliptic case one finds that
\begin{equation} \ell = \sqrt{x_o(R_A - s_o)} = \sqrt{(L -s_o)
(R_B - L + s_o)}. \end{equation}
These conditions uniquely determine $x_o, \ell$ (see the pictures on p.135 of [L.1]).
It also follows that $c_{00}(s) = - \frac{\psi_{11}}{\psi_{11}} = - \frac{\frac{s-s_0}{\ell}}
{\pm \ell+ \frac{(s-s_0)^2}{\ell} }.$

In \S 4.2 we will put the leading term into a canonical normal form with coefficients
independent of $\s.$ 

\subsection{Wave invariants revisited}

Having straightened the Laplacian, and hence having transferred the information about
the boundary into the metric, 
 we can now state precisely just how much data of the
boundary goes into a given wave invariant $a_{\gamma k}.$   Since the wave invariants
at $\{\y = 0\}$ of $\bDelta$ are calculated by the stationary phase method, we have,
as in the boundaryless case:
\begin{prop} \label{BDROP} 
$a_{\gamma k}(\bDelta) = a_{\gamma k}(\bDelta_1^{\leq 2(k+2)}  + \bDelta_o^{\leq 2(k+1)}
 + \dots + \bDelta_{-k -1}^{o}).$ \end{prop}

As $\bDelta$ is a well-defined partial differential operator in a neighborhood of
$\Omega_o$ we can again reformulate the conclusion in terms of the symbol classes
$S^{2, k}(V, \R^+\gamma_o)$ where  $\gamma_o$ denotes the bouncing ball orbit $[0, L] \times \{0\}$.
The discussion in the boundaryless case remains valid, althought the microlocal
neighborhood $V$ now acquires two components: 
Since $|\sigma| \geq \half > 0$ in  $V$,   $V = V_+ \cup V_-$
where $V_+ = V \cap \{\sigma > 0\}, V_- = V \cap \{\sigma < 0\}.$ The two components
are obviously interchanged by the canonical involution $\tau (s, \sigma, y, \eta) =
(s, - \sigma, y, -\eta)$ of $T^*([0, L] \times \R)$.

 \begin{cor}\label{CORBDROP} 
$a_{\gamma k}(\bDelta) $ depends only on the class of $\bDelta$ modulo $S^{2, 2(k + 2)}$. 
 \end{cor}

In view of the metric normal form,  the data which goes in to the k-jet of $\bDelta$ along
$\{y = 0\}$ is precisely the $k+ 2$-jet of the function $B(\s, \y)$ and hence the $k + 2$-jet of the boundary defining functions.  Hence 
$a_{\gamma k}$ depends only on the $2k + 4$-jet of the boundary at $\{y = 0\}.$

We will actually need a slight generalization of (\ref{CORBDROP}) which is prove in precisely the same way.

\begin{prop}\label{FIODROP} Suppose that $A \in Op (S^{2, 2(k + 2)})$.  Then $res F A G = 0$ for any bounded Fourier integral
operators $F, G.$ \end{prop}

\section{Semiclassical normal form}

As mentioned in the introduction, our approach is to convert the conjugation to normal form to a semiclassical
problem.   

\subsection{Semiclassical scaling}

To  introduce the semiclassical parameter, we  make a  semiclassical scaling of operators on the model space.
  Roughly, the scaling weights the
tangential derivative $D_{\s}$ by $N^2$, the normal derivative $D_{\y}$ by $N$ and
$\y$ by  $N^{-1}.$

 We define
operators $T_N, M_N$ on the model space $L^2(\Omega_o)$ by
\begin{equation} \begin{array}{ll} \bullet &  T_N f(\s,\y):= N f(\s, N \y) \\
\\ \bullet &  M_N f(\s,\y) := e^{i N^2 \s} f(\s,\y) \end{array} \end{equation} We then
have: \begin{equation} \begin{array}{l}  T_N^* D_{\y} T_N = N  D_{ \y} \\ \\
T_N^* \y  T_N  =
N^{-1} \y \\ \\
 M_N^*D_{\s} M_N =( N^2 + D_{\s}) \end{array} \end{equation}

\begin{defn} The rescaling of an operator $a^w(\s,D_{\s},\y,D_{\y})$ is given by 
$$a^w_N(\s,D_{\s},\y,D_{\y}) :=M_N^* T_N^*a^w(\s,D_{\s},\y,D_{\y})T_N M_N.$$ 
\end{defn}

We have:
\begin{prop} The complete symbol $a^w_N(\s, \bsigma, \y, \bareta)$ of $M_N^* T_N^*a^w(\s,D_{\s},\y,D_{\y})T_N M_N$ is given by
$$a^w_N(\s,\bsigma,\y,\bareta) = a^w(\s, \bsigma +  N^2, \frac{1}{N} \y, N \bareta).$$
\end{prop}

\noindent{\bf Proof} The operator kernel of $a^w(\s, D_{\s}, \y, D_{\y})$ is equal to
$$\int_{\R}\int_{\R} a(\half(\s - \s'), \bsigma, \half(\y - \y'), \bareta) e^{i ( \langle \s - \s', \bsigma \rangle + \langle \y - \y', \bareta\rangle }
d \bsigma d\bareta.$$
Conjugating with  $M_N$ amounts to adding $N^2 (\s - \s')$ to the phase and hence to a translation $\sigma \to \sigma + N^2$. Conjugating
with $T_N$ amounts to changing the phase $ \langle \y - \y', \bareta\rangle$  to $\frac{1}{N}  \langle \y - \y', \bareta\rangle$
and the amplitude to $a(\half(\s - \s'), \bsigma, \half(\frac{\y}{N} - \frac{\y'}{N}), \bareta).$ Change variables $\bareta \to N \bareta$ to get a Weyl pseudodifferential
operator with amplitude $a(\half(\s - \s'), \bsigma, \half(\frac{\y}{N} - \frac{\y'}{N}), N \bareta).$ \qed

\begin{defn} The semiclassically scaled (1/2-density) Laplacian is the Weyl pseudodifferential operator
on $\Omega_o$ defined by
$$\bDelta_N = M_N^* T_N^* \bar{\Delta} T_N M_N .$$ \end{defn}

In the straightened form, we have
\begin{equation}\begin{array}{l} \bDelta_N \sim (D_{\s} + N^2)^2  + 
B(s, N^{-1} \y) [\pm  \y^2 N^{-2} (D_{\s} + N^2)^2 + N^2 D_{\y}^2] \\ \\ + \{\Gamma_{\s}\}_N  [D_{\s} + N] + N \{\Gamma_{\y}\}_N  D_{\y}
+ \{K\}_N \end{array} \end{equation}
where $\{f(\s, \y)\}_N = T_N^* f(s,y) = f(s, N^{-1}y)$ and where (as usual) the alternative $\pm$ refers to the
elliptic/hyperbolic dichotomy. 

The Weyl symbol of $\bar{\Delta}$ has the simple form:
\begin{equation} \sigma_{\bar{\Delta}}^w(\s, \bar{\sigma}, \y, \bar{\eta}) :=
\sigma^2 + B(s, y) I + K,\;\;\;\; K = \bar{\Delta}\cdot 1\end{equation}
The linear terms vanish because they give the subprincipal symbol in the
Weyl calculus and that of $\Delta$ equals zero.  By the above proposition,
we get upon rescaling
\begin{equation}\sigma_{\bar{\Delta}_N}^w   
 \sim   (\sigma + N)^2 + B(s, \frac{1}{N} y) (\pm y^2 N^{-2}(\sigma + N^2)^2 + N^2 \eta^2)
+ K(s, \frac{1}{N} y).  \end{equation}

\subsection{Linearized problem}

To get a sense of what is involved in putting $\Delta$ into normal form, let
us first consider the `linearized' problem which involves only the highest
powers of $N$: 
\begin{equation} \sigma_{\bar{\Delta}_N}^w  =  N^4 + 2 N^2 [  D_{\s} 
+    \dot{b}_{00}(s) \hat{I}] \;\;\rm{mod}\;\;N
\end{equation}
where $\hat{I}$ denotes the quantum action operator:  $\hat{I}^e = \half (D_{\y}^2 +  y^2)$ in the
elliptic case and  $\hat{I}^h = \half (D_{\y}^2 -  y^2)$ in the hyperbolic case. 

As in
 the boundaryless case [Z.1, \S 1] we can complete
the conjugation of the linear/quadratic term to normal form by a moving metaplectic conjugation. When
results apply mutatis-mutandi to both elliptic and hyperbolic cases, we denote the action simply by $\hat{I}.$
In the following proposition, $\mu$ denotes the metaplectic representation of $SL(2, \R)$ (strictly speaking,
of its double cover but signs are irrelevant here).

\begin{prop} There exists an $SL(2, \R)$-valued function $a_{\alpha}(\s)$ so that 
$$\begin{array}{l} \mu(a_{\alpha})^* [ D_{\s} + b_{00}(\s) \hat{I}] \mu(a_{\alpha}) =
{\cal R}\\ \\\mu(a_{\alpha})(0) = \mu(a_{\alpha}) (L) = Id \end{array}$$  where 
${\cal R} =  D_s + \frac{\alpha}{L} \hat{I}$ and where $\alpha = \int_0^L b_{00}(s) ds.$
\end{prop}

\noindent{\bf Proof} We first construct a function $a(s)$ so that $ \mu(a)^* [ D_{\s} + b_{00}(\s) \hat{I}] \mu(a) =
D_{\s}.$  
 The desired metaplectic operator is obviously given by  
\begin{equation} \mu(a(\s)) = \exp (-  i [\int_0^{\s} b_{00}(\s) d\s] \hat{I}).
\end{equation}
Hence $a(\s) =  \exp (- J/ 2  \int_0^{\s} b_{00}(\s)  d\s)$ where 
$$\begin{array}{ll} J = \left( \begin{array}{ll} 0 & 1 \\ -1 & 0 \end{array} \right) & \rm{elliptic\; case} \\  &\\
J = \left( \begin{array}{ll} 1 & 0 \\ 0 & -1 \end{array} \right) & \rm{hyperbolic\; case}.\end{array} $$ 

  The boundary values  of $a$ are
given by $a(0) = I,\;\;\;\; a(L) = \exp (- J
  \int_0^L b_{00}(\s)  d\s) = \exp(- \alpha I)$.
Now let $r_{\alpha}(s) = \exp ( \frac{\alpha}{L} J \s) $ and put $a_{\alpha}(s) = a(s) r_{\alpha}(s)^{-1}.$  Then  $\mu(a_{\alpha}(0)) = Id$ 
and $\mu(a){\alpha}(L) =  \exp( \alpha \hat{I}) \circ \exp(- \alpha \hat{I}) = Id. $ Further, $\mu(r_{\alpha}(s))$
conjugates $D_s$ to $D_s + \frac{\alpha}{L} \hat{I}.$ \qed

\subsection{ Semiclassical pseudodifferential operators}\label{SCO}

 Semiclassical scaling produces a partial differential operator $\bDelta_N$ depending
on a small parameter $1/N$ and in the conjugation to normal form we will introduce other such operators.  Let us pause to clarify the kinds of semiclassical pseudodifferential operators
which will be of concern to us.  Our discussion is based on ideas and notation from \cite{BMGH}.

First, let us recall that 
the usual (admissible) semiclassical  $h = 1/N$-pseudodifferential operators of order m on $\R^n$ are the Weyl
quantizations
$$a^w(x, \frac{1}{N} D; N) u(x) =(2 \pi h)^{-n} \int\int e^{i \langle x-y, \xi \rangle}
a(\frac{1}{2}(x+y), \frac{1}{N} \xi; N) u(y) dy d\xi $$
of  amplitudes $a$ belonging to the space $S_{cl}^m(T^*\R^n)$ of asymptotic
sums
$$a(x, \xi, N) \sim N^{ m}\sum_{j =0 }^{\infty} a_j(x, \xi) N^{- j}   $$
with $a_j \in C^{\infty}(R^{2n}).$ 

   Semiclassical scaling gives rise to symbols of the  form
$a(\s, \bsigma +  N^2, \frac{1}{N} \y, N \bareta)$ where $a $ is a polyhomogeneous
symbol $a \sim \sum_{j = 0}^{\infty} a_{m - j}$  of order $m$. 
As in (\ref{HOMO}) we may write:
\begin{equation} \label{SCHOMO}  a(\s, \sigma  + N^2 ,  \frac{1}{N} \y, N \bareta)
 \sim \sum_{j = 0}^{\infty} N^{2m - 2j} a_{m - j}(\s, 1 + \frac{\sigma}{N^2}, \frac{1}{N} \y, \frac{1}{N} \bareta).  \end{equation} 
It is evident that semiclassical scaling produces symbols which behave  in both transverse
variables $(\y, \bareta)$ like standard semiclassical symbols in the fiber variable
$\eta$.  Thus, scaled symbols are isotropic analogues of semiclassical admissible symbols
in the transverse $\R^n$. They also have a tangential dependence in $1 + \frac{\sigma}{N}$
which seems to have  no precise analogue for standard semiclassical symbols. 

 Let $\psi(\s, \sigma, \y, \bareta)$ be the homogeneous cutoff to the cone $V$.  Under
rescaling it goes over to the symbol $N^2 \psi(\s, 1 + \frac{\sigma}{N^2}, \frac{\y}{N},
\frac{\bareta}{N})$.  For symbols independent of $\sigma$ or for $\sigma < \epsilon N^2$
the scaled symbol is supported in the transverse ball
\begin{equation} V_N:= \{(\s, \sigma, \y, \bareta): |(\frac{1}{N} \y, \frac{1}{N} \bareta) | \leq \epsilon \}. \end{equation}

It is useful to reformulate the condition that $A \in Op S^{m, k}$ in terms of the scaled symbol.
  Given $A \in Op S^m$ we define (with some modifications to \cite{BMGH})  the formal
differential operator
\begin{equation} \sigma^{\infty}_{(\s, \sigma, x, \xi)} (A) = \sum_{j, \alpha, \beta} \frac{1}{\alpha ! \beta!} (\frac{\partial}{\partial x})^{\alpha} 
(\frac{\partial}{\partial \xi})^{\beta} a_j(\s, \sigma, x,\xi) \y^{\alpha} D_{\y}^{\beta} N^{- (|\alpha| + |\beta| + 2j)} \end{equation}
and put
\begin{equation} \sigma^{k}_{(\s, \sigma, x, \xi)}(A) = \sum_{|\alpha| + |\beta| + 2j = k} \frac{1}{\alpha ! \beta!} (\frac{\partial}{\partial x})^{\alpha} 
(\frac{\partial}{\partial \xi})^{\beta} a_j(\s, \sigma, x, \xi) \y^{\alpha} D_{\y}^{\beta}.\end{equation}
Then $a \in S^{m, k}(V, \R^+\gamma)$ if and only if 
$\sigma_{(\s, \sigma, x, \xi)}^{\infty}$ is divisible by $N^{-k}$ for  $(x, \xi) = (0,0)$; hence $\sigma_{(\s, \sigma, \y, \bareta)}^{\infty}(A) = \sigma^k_{(\s, \sigma, \y, \bareta)}(A) N^{-k}$ mod $N^{-(k + 1)}$ for
$(s, \sigma, x, \xi) \in \R^+ \gamma_o$.  
 We thus have
\begin{prop}\label{SYM} If $a \in S^m(V)$ then  $a \in S^{m, k}(V, \R^+ \gamma_o)$ if and only if the formal Taylor expansion of   $a(\s, N^2, \frac{1}{N} \y, N \bareta)$ along $\R^+ \gamma_o$ is
divisible by $N^{-(k + 1)}.$ \end{prop}

\subsection{Conjugation to a semiclassical normal form}

We now come to the principal step in the conjugation to normal form: the conjugation
of $\bDelta_N$ to a semiclassical normal form.  

We first conjugate by $\mu(a_{\alpha})$ as in
 the linearization step to get the somewhat simpler form:
\begin{equation} {\cal R}_N : = \mu(a_{\alpha})^* \bDelta_N \mu(a_{\alpha}). \end{equation}
We then wish to  conjugate  ${\cal R}_N$ to the semiclassical normal form 
\begin{equation} \label{SCNF} F_N(\hat{I})^2 \sim N^4 + N^2 \frac{\alpha \hat{I}}{L} + p_1(\hat{I}) + 
N^{-2} p_2(\hat{I}) + \cdots \end{equation}
 by means of a semiclassical
pseudodifferential intertwining operator $W_N(\s, \y, D_{\y})$ which preserves Dirichlet boundary conditions. 
Thus the full intertwining operator is $  \mu(a_{\alpha}) W_N.$ To avoid encumbring
the notation we will also denote this full intertwining operator in  (\S 5)by $W_N$.

Let us now explain what we mean by semiclassical conjugation to normal form. We assume for simplicity that the
bouncing ball orbit is elliptic, but the same  argument and result hold in the hyperbolic case.
 Roughly speaking, our  object is to produce  bounded semiclassical
pseudodifferential operators $\tilde{W}_k^+(\s, \y, D_{\y}), \tilde{W}_k^-(\s, \y, D_{\y})$ defined in the open neighborhood (even $U_{\epsilon}$ of $\Omega_o$ and satisfying the following asymptotic relations on this domain:
\begin{equation} \label{QM} \begin{array}{ll} (i) & \bDelta \tilde{W}_k^+ e^{i \pi k \frac{\s}{L}} D_q(N_k \y) \sim F(k, q + \half)^2 e^{i \pi k \frac{\s}{L}} \tilde{W}_k^+  D_q(N_k \y) \\ & \\
(ii) & \bDelta \tilde{W}_k^- e^{- i \pi k \frac{\s}{L}} D_q(N_k \y) \sim F(k, q + \half)^2 e^{- i \pi k \frac{\s}{L}} \tilde{W}_k^- D_q(N_k \y) \\ & \\
(iii) & \tilde{W}_k^+ e^{i \pi k \frac{\s}{L}} D_q(N_k \y) - \tilde{W}_k^- e^{- i \pi k \frac{\s}{L}} D_q(N_k \y) = 0\;\; \mbox{at}\;\;
\s = 0, \s = L. \end{array} \end{equation}
Here, $F(k, q+ \half) = F_{N_k}(q + \half)$ with $N_k = \sqrt{k}$ and $D_q$ is the qth  normalized Hermite function. The precise meaning of $\sim$ will be clarified
below.  The Laplacian $\bDelta$ is the Laplacian acting on the open space $U_{\epsilon}.$  The condition (iii) implies
that $\tilde{W}_k^+ e^{i \pi k \frac{\s}{L}} D_q(N_k \y) - \tilde{W}_k^- e^{- i \pi k \frac{\s}{L}} D_q(N_k \y) \in H^1_0(\Omega)$
so that it lies in the domain of $\bDelta_{\Omega}$.

Since $C \bDelta C = \bDelta$ (with $C$ the operator of complex conjugation), it suffices to construct $\tilde{W}^+_k$
and to 
put $\tilde{W}_k^- = C \tilde{W}_k^+ C$.  Then, 
\begin{equation} \tilde{W}_k^+ e^{i \pi k \frac{\s}{L}} D_q(N_k \y) -
 \tilde{W}_k^- e^{- i \pi k \frac{\s}{L}} D_q(N_k \y) = 2i \Im \tilde{W}_k^+ e^{i \pi k \frac{\s}{L}} D_q(N_k \y) \end{equation}
In order that $\Im \tilde{W}_k^+ e^{i \pi k \frac{\s}{L}} D_q(N_k \y)$ lie ${\cal D}om (\bDelta_{\Omega})$ it is thus
sufficient that  
\begin{equation}\label{C} C \tilde{W}_k^+(0) C = \tilde{W}_k^+(0), \;\;\; C \tilde{W}_k^+(L) C = \tilde{W}_k^+(L). \end{equation}
This boundary condition (\ref{C}) is correct in both the elliptic and hyperbolic cases,
although the  quasimode construction which motivates it  only works in the elliptic case.  

As a further preliminary, let us  rewrite the equations (\ref{QM}) and the boundary conditions (\ref{C}) in terms
of the semiclassically scaled Laplacian.  We observe that $D_q(N_k \y) = T_{N_k} D_q(\y)$. Let us define:
\begin{equation}\label{W} W_k^+ = T_{N_k}^{-1} \tilde{W}_k^+ T_{N_k}. \end{equation}  Then (\ref{QM})-(\ref{C}) is equivalent to
 \begin{equation} \label{QMS} \begin{array}{ll} (i) & \bDelta_{N_k} W_k^+  D_q( \y) \sim F(k, q + \half)^2 W_k^+ D_q( \y) \\ & \\
(ii) & \bDelta_{N_k} W_k^-  D_q( \y) \sim F(k, q + \half)^2 W_k^- D_q( \y) \\ & \\
(iii) & C W_k^+(0) C = W_k^+(0), \;\;\; C W_k^+(L) C = W_k^+(L). \end{array} \end{equation}

The operators   $W_k^+ = W_k^+(\s, \y, D_{\y})$ will be  essentially a family of pseudodifferential
operators on the transverse space, parametrized by $\s.$  Hence there are no subtleties involving the definition
of pseudodifferential operators on manifolds with boundary.  Moreover,  $W_k^+$ is essentially applied
only to a function of $\y.$ To be more precise, we first conjugate by $\mu(a_{\alpha})$ to put the linear term
in normal form, and then $W_k^+$ is applied to a function in the kernel of ${\cal R}$. 
Hence we only require that the conjugation identity hold as operators applied to
 functions in the kernel of ${\cal R}$. 
 Therefore we introduce the following notation:
   Given an operator $A(s, D_s, y, D_y)$, we denote by $A|_o$ the restriction of $A$ to 
functions of $y$ only, i.e. $|_o$ denotes the restriction to functions in the kernel of ${\cal R}$.

The following lemma proves the existence of such a conjugating operator.  We emphasize that the conjugation of
$\bDelta$ takes place over $U_{\epsilon}$.

\begin{lem} \label{MAIN} Let ${\cal R}_N^+ = \mu(a_{\alpha})^* \bDelta_N \mu(a_{\alpha})$ with 
  $\frac{\alpha}{\pi} \notin \Q$ in the elliptic case.  Then there exist  polynomial differential  operators $ P_{j/2}^w
(\s, \y, D_{\y})$ and $Q_{j/2}^w(s, \y, D_{\y})$ of 
degree $2j + 2$ 
on $L^2(\R_y)$ with smooth coefficients in $\s \in (-\epsilon, L + \epsilon)$ and polynomials $f_j(\hat{I})$ of
degree $j+2 $   such that:

\noindent(a)  $ P_{j/2}(\s, \y, \eta)$ and $ Q_{j/2}(\s, \y, \eta)$  are real-valued;

\noindent(b) For each $\s$, there  is a formal $N$-expansion: $$(W_{N \frac{K}{2}}^{+ })^{-1} {\cal R}_N  W_{N \frac{K}{2}}^{+} \sim - N^4  + 2 N^2 {\cal R} +
\sum_{j=0}^{\infty} N^{-j}
 {\cal R}^{\frac{K}{2}}_{2-\frac{j}{2}}(\s,D_{\s},\y,D_{\y})$$
where
\medskip

(i) $ {\cal R}^{\frac{K}{2}}_{2-\frac{j}{2}}(\s, D_{\s}, \y, D_{\y}) =
  {\cal R}^{\infty,2}_{2-\frac{j}{2}} {\cal R}^2 +
   {\cal R}^{\infty,1}_{2-\frac{j}{2}} {\cal R} + {\cal R}^{\infty,o}_{2-\frac{j}{2}},$ with
${\cal R}^{\infty,k}_{2-\frac{j}{2}} \in C^{\infty}([0, L], {\cal E}_{\epsilon}^{j -2k});$
\medskip

(ii) $ {\cal R}^{\frac{K}{2}}_{2-j}(\s, D_{\s}, \y, D_{\y})|_o = {\cal
R}^{\infty,o}_{2-j}(\s, \y, D_{\y})|_o =
  f_j(\hat{I})$ for certain polynomials $f_j$
 of degree j+2 on $\R$;
\medskip

(iii) $ {\cal R}^{\frac{K}{2}}_{2-\frac{2k+1}{2}}
(\s, D_{\s}, \y, D_{\y})|_o = {\cal R}^{\infty, o}_{2-\frac{2k+1}{2}}(\s,\y,
D_{\y})|_o = 0;$
\medskip

\noindent(c) $W_{N \frac{K}{2}}^{+ }$ satisfies the boundary condition: $C W_{N \frac{K}{2}}^+(0) C = W_{N \frac{K}{2}}^+(0), \;\;\;\;\;\; CW_N^+(L) C = W_{N \frac{K}{2}}^+(L); $

\noindent(d) We have $$[(W_{N  \frac{K}{2}}^{+})^{-1} {\cal R}_N  W_{N \frac{K}{2}}^{+}]|_o \sim - N^4  + 2 N^2 {\cal R} + F_{N K}^{2 }(\hat{I})+ E_{N \frac{K}{2}}$$  where
$F_{N K}^{ 2}(\hat{I}) : = \sum_{j=0}^{K} N^{-2j} f_j(\hat{I})$ and 
where 

\noindent(e) The error term $E_{N +\frac{K}{2}}(\s, \y,  \bareta) = N^{-(K+1)} \tilde{E}_{N \frac{K}{2}}$ where $\hat{I}^{-(K + 1)}\tilde{E}_{N \frac{K}{2}}$
is a bounded operator on $L^2(\R)$ for each $\s$ with a uniform bound in $N$.    

\end{lem}

\noindent{\bf Proof}  

After the linearization step, the scaled Laplacian has the perturbative form
\begin{equation} {\cal R}_N^+ = N^4 + N^2 {\cal R} + E_{N 0}^{+},\;\;\;\; E_{N 0}^{+} := \sum_{m=2}^{\infty} N^{(4 - m)}
{\cal R}_{2 -\frac{m}{2}} \end{equation} 
where
\begin{equation} {\cal R}_{2 -\frac{m}{2}} \in C^{\infty}(\R,{\cal E}^{m-4}_{\epsilon}){\cal R}^2 
 +C^{\infty}(\R,{\cal E}^{m-2}_{\epsilon}) {\cal R}
+C^{\infty}(\R,{\cal E}^m_{\epsilon})\end{equation}
Actually, as discussed in (\ref{DROP}), $a_{\gamma k}$ depends only on  the class of $\bDelta$
in $S^{2, 2(k + 2)}(T^*U, \R^+ \gamma)$ and by (\ref{BDROP}) this is the same as the class
of the complete symbol $\sigma_{\bDelta_N}$ of $\bDelta_N$ at $\sigma = 0$ modulo
$N^{-2(k + 2)}.$  We may therefore drop higher order terms 
 to get a polynomial partial differential operator $\bDelta^{2k + 2}$ 
with the same wave invariants $a_{\gamma j}$ for $j \leq k.$  . For simplicity we do not indicate this
truncation in our notation, but the reader is invited to think of $\bDelta$ as a polynomial differential operator
in the $(\y, D_{\y})$ variables. 

We now  construct Weyl symbols   $P_{j/2},  Q_{j/2}$ so that iterated composition with $e^{ N^{-j}(P + i Q)_{j/2}^w}$ will successively remove the lower order terms in ${\cal R}_N$  after  restriction by $|_o$ and so that the
boundary condition is satisfied.
As mentioned above, all operators will be standard Weyl pseudodifferential operators
defined in a neighborhood of $\Omega_0$ and  acting only on the transverse space $\R_y$ with coefficients in $\s$.
Hence we may construct them using the  symbolic calculus.  

Let us first rewrite  the boundary condition in terms of Weyl symbols:
On the symbol level, we have
\begin{equation} C a^w(\s, \y, D_{\y}) C = \bar{a}^w (\s, \y, - D_{\y}) \end{equation}
where $\bar{a}$ is the complex conjugate of $a.$  Hence it is natural to split up a Weyl symbol into its real/imaginary
parts and into its even/odd parts with respect to the canonical involution $(\y, \bareta) \to (\y, - \bareta)$.
We note that this splitting is invariantly defined since $a$ is real if and only if $a^w(\y, D_{\y})$ is self-adjoint
and since the even/odd parts of a real symbol are its even/odd parts under conjugation by $C$. Since our symbols
are always Weyl symbols in this section, we omit the superscript $w$ in the future.

By assumption $P_{j/2}$ and $Q_{j/2}$ are real symbols so it remains to split them into their even/odd parts:
\begin{equation} P_{j/2} = P_{j/2}^e + P_{j/2}^o,\;\;\;\; Q_{j/2} = Q_{j/2}^e + Q_{j/2}^o. \end{equation}
Now we observe that for general Weyl pseudodifferential operators $P, Q$,  
$$ C (P + i Q) C =  i Q \;\;\mbox{iff}\;\; \; P^o = Q^{e} = 0. $$
Therefore, the boundary condition on $W_N^+$ is equivalent to:
\begin{equation} \label{BC} \begin{array}{l}  P_{j/2}^o(0, \y, \bareta) = P_{j/2}^o(L, \y, \bareta) = 0,\;\;\;\;\;\;
 Q_{j/2}^e(0, \y, \bareta) = Q_{j/2}^e(L, \y, \bareta) = 0. \end{array} \end{equation}
We emphasize that there is no condition on $Q_{j/2}^o$ or $P^e_{j/2}$.

To begin the induction, let us construct  $P_{\frac{1}{2}}(\s, \y, D_{\y}),\; Q_{\frac{1}{2}}(\s, \y, D_{\y}) \in C^{\infty}([0,L])\otimes{\cal E}^3_{\epsilon}$ so that the boundary conditions are satisfied and so that
\begin{equation} e^{-  N^{-1}(P + i  Q)_{\frac{1}{2}}} {\cal R}_N e^{ N^{-1} (P+ i Q)_{\frac{1}{2}}}|_o = \{ N^4 + N^2 {\cal R} + E_N^1 \}|_o. \end{equation} 
Expanding the exponential, we get to leading order the {\it homological equation}: 
\begin{equation}  \{ [{\cal R} , (P+ i Q)_{\frac{1}{2}}]+ {\cal R}_{\frac{1}{2}} \}|_o= 0. \end{equation}
Taking the complete symbol of both sides we get the symbolic homological equation:
\begin{equation} i \{\sigma + \frac{\alpha}{L} I, P_{\half} +  i Q_{\frac{1}{2}}\}+ {\cal R}_{\frac{1}{2}}|_o = 0.\end{equation}  
The equation may be rewritten in the form:
\begin{equation} \label{SOL} \partial_{\bar{s}} (P + i Q)_{\frac{1}{2}}(\s, r_{\alpha}(\s)(\y, \bareta))=
 -i  {\cal R}_{\frac{1}{2}}|_o.   \end{equation} 
The solution has the form:
\begin{equation} (P + i Q)_{\half}(\s,r_{\alpha}(\s)(\y, \bareta))  =   (P + i Q)_{\half}(0) - 
 i \int_0^{\s}  {\cal R}_{\half}(u, \y, \bareta) du.\end{equation}
We need to determine $ (P + i Q)_{\half}(0)$ so that the boundary conditions at $\s = 0$ and $\s = L$ are
satisfied.  To clarify the boundary conditions and their solvability, we
separate out  real/imaginary and even/odd parts in the equations to get:
\begin{equation} \label{EO} \begin{array}{ll} \partial_{\s} P_{\half}^e + \frac{\alpha}{L}  \{I, P^o_{\half}\} =  \Im {\cal R}_{\frac{1}{2}}^e|_o &
\partial_{\s} P_{\half}^o + \frac{\alpha}{L}  \{I, P^e_{\half}\} =  \Im {\cal R}_{\frac{1}{2}}^o|_o 
\\ & \\
\partial_{\s} Q_{\half}^e + \frac{\alpha}{L}  \{I, Q^o_{\half}\} = -\Re  {\cal R}_{\frac{1}{2}}^e|_o &
\partial_{\s} Q_{\half}^o + \frac{\alpha}{L}  \{I, Q^e_{\half}\} = - \Re {\cal R}_{\frac{1}{2}}^o|_o \end{array} \end{equation}
together with the boundary conditions
\begin{equation} P_{\half}^o(0) = P^o_{\half}(L) = 0,\;\;\;\;\;\;Q_{\half}^e(0) = Q_{\half}^e(L) = 0.\end{equation}

 We now observe as in the boundaryless case (cf. \cite{Z.1}, Lemma (2.22))  that ${\cal R}_{\half}^{\pm}(u, \y, \bareta ) = \y \circ (a r_{\alpha}^{-1} (u)) I$ is a   polynomial of degree 3 in $(\y, \eta)$ in which every term is of odd degree in $(\y, \eta)$.   By (\ref{SOL}) it follows that $P_{\half}, Q_{\half}$ are  also odd polynomials
of degree 3.

To analyse the equations further we change coordinates. In the elliptic case, we use 
 the complex cotangent variables $z = \y + i \bareta,
\bz = \y - i \bareta$ which satisfy  $\{I, z^m \bz^n\} = i ( z \frac{\partial}{\partial z} - \bz \frac{\partial}{\partial \bz}) z^m \bz^n =
i (m - n) z^m \bz^n$. 
We then write:
\begin{equation}\begin{array}{l}P_{\half}^e(\s, z, \bz) = \sum_{m,n: m + n \leq 3} p_{\half mn}^e(\s) (z^m \bz^n +
\bz^m z^n), \\ \\
P_{\half}^o(\s, z, \bz) = \sum_{m,n: m + n \leq 3} p_{\half mn}^o(\s) (z^m \bz^n -
\bz^m z^n). 
\\ \\  
Q_{\half}^e(\s, z, \bz) = \sum_{m,n: m + n \leq 3} q_{\half mn}^e(\s) (z^m \bz^n +
\bz^m z^n), \\ \\
Q_{\half}^o(\s, z, \bz) = \sum_{m,n: m + n \leq 3} q_{\half mn}^o(\s) (z^m \bz^n -
\bz^m z^n). 
\end{array} \end{equation}
In the hyperbolic case, we first  use real coordinates  $(y, \eta)$ in which $I^h = \half( \eta^2 - y^2)$.
We then define $w = y + \eta, \bar{w} = y - \eta$ so that $I^h = \half |w|^2 := \half w \bar{w}$ 
and so that   $\{I^h, w^m \bar{w}^n\} = (m - n) y^m \eta^n$.  Since the elliptic and hyperbolic cases are quite
similar, we 
only carry out the details in the elliptic case and refer to \cite{Z.2} (Lemma 3.1) for complete details on
the hyperbolic case. We will refer to $m = n$ terms  as the `diagonal terms'.  They only occur for even $j$ and are terms which
are polynomials in the action variables.

Then (\ref{EO}) may be rewritten in terms of these coordinates. In the elliptic case, we
have:
\begin{equation} \label{EOC} \begin{array}{l}
\frac{d}{ds} p^e_{\half mn}(s) +  i \frac{\alpha}{L}  (m - n) p^o_{\half mn}(s)  =   \Im {\cal R}_{\frac{1}{2}}^e|_o \\ \\
 \frac{d}{ds} p^o_{\half mn}(s) + i \frac{\alpha}{L} (m - n)  p^e_{\half mn}(s)  =  \Im  {\cal R}_{\frac{1}{2}}^o|_o \\ \\
 \frac{d}{ds} q^e_{\half mn}(s) +  i \frac{\alpha}{L}  (m - n) q^o_{\half mn}(s)  = - \Re {\cal R}_{\frac{1}{2}}^e|_o \\ \\
 \frac{d}{ds} q^o_{\half mn}(s) + i \frac{\alpha}{L} (m - n)   q^e_{\half mn}(s)  = - \Re {\cal R}_{\frac{1}{2}}^o|_o \end{array} \end{equation}
Similarly in the hyperbolic case, although there is no factor of $i$ in the second terms on the left sides. 
Since $m \not= n$, we  can immediately solve for $p^e_{\half mn},  q^o_{\half mn}$:
\begin{equation} \begin{array}{l} p^e_{\half mn}(s)  = \frac{- L}{i  \alpha (m - n)} \{\frac{d}{ds} p^o_{\half mn}(s)  -
 \Im  {\cal R}_{\frac{1}{2}}^e|_o\}\\ \\
 q^o_{\half mn}(s)  = \frac{- L}{i  \alpha (m - n)} \{\frac{d}{ds} q^e_{\half mn}(s)  +
  \Re {\cal R}_{\frac{1}{2}}^e|_o\}
\end{array} \end{equation}  
We thus eliminate the $p^e_{\half mn}, q^o_{\half mn}$ variables 
and reduce to  (uncoupled) second order equations for the independent variable $p^o_{\half mn}, q^e_{\half mn}(s).$
In the elliptic case, they read:
\begin{equation} \label{2O} \begin{array}{l} - \frac{d^2}{d s^2}p^{o}_{\half mn}(s) -  [\frac{\alpha}{L} (m - n)]^2 p^{o}_{\half mn}(s) =
 \frac{L}{\alpha (m-n)}  \frac{d}{ds} \Im  {\cal R}_{\frac{1}{2}}^e|_o +   \Im {\cal R}_{\frac{1}{2}}^e|_o
\\ \\ 
- \frac{d^2}{d s^2}q^{e}_{\half mn}(s) -  [\frac{\alpha}{L} (m - n)]^2 q^{e}_{\half mn}(s) =
- \frac{L}{\alpha (m-n)}  \frac{d}{ds}  \Re {\cal R}_{\frac{1}{2}}^e|_o +  \Re  {\cal R}_{\frac{1}{2}}^e|_o
 \end{array} \end{equation}
In the hyperbolic case, the operator on the left side is $- \frac{d^2}{ds^2} + [\frac{\alpha}{L} (m - n)]^2.$ 
The boundary conditions on $p^o_{\half mn}, q^e_{\half mn}$ are (in both elliptic and hyperbolic cases):
\begin{equation}\label{bc} \begin{array}{ll}  p^o_{\half mn}(0)  =  0, & p^o_{\half mn}(L)  =  0 \\ & \\  
q^e_{\half mn}(0)  =  0, & q^e_{\half mn}(L)  =  0. \end{array} \end{equation}
The  boundary value problem (\ref{2O}) - (\ref{bc}) is always solvable  unless $0$ is an eigenvalue of the operator
$ D_{\s}^2 - [\frac{\alpha}{L}  (m - n)]^2$ (elliptic case), resp. $ D_{\s}^2 + [\frac{\alpha}{L}  (m - n)]^2$ (hyperbolic case) with boundary conditions $q(0) = 0 = q(L)$  The eigenfunction
would have to have the form $\sin (\frac{\alpha}{L}  (m - n) \s)$ in the elliptic case or $\ sinh (\frac{\alpha}{L}  (m - n) \s)$ in the hyperbolic case.  Hence in the elliptic case, a sufficient condition for solvability is that $\alpha/\pi \notin  \Q$ while in the hyperbolic case there is no obstruction if $\alpha \not= 0.$ 

Thus we have solved the conjugation problem to third order.  The exponents $P_{\half}, Q_{\half}$   
are odd polynomial differential
operators of degree 3  with smooth  coefficients defined in a neighborhood of $[0, L]$.  By construction, the $Q_{j/2}$'s  will always have  the same order,  same order of vanishing, and same parity 
as the restriction  ${\cal R}_
{\frac{1}{2}}|_o$.  It follows that 
$$N^{-1}  ad((P + i Q)_{\half}): N^{- k}\Psi^l(\R) \otimes {\cal E}^m_{\epsilon}
\rightarrow N^{-(k+1)} [\Psi^{l-1}(\R) \otimes {\cal E}^{m+3}_{\epsilon} + 
\Psi^{l}(\R)
\otimes {\cal E}^{m+1}_{\epsilon}].$$

We now carry the process forward one more step because, as in the boundaryless case, the
even steps require something new.
In  the second step, ${\cal R}_N$ is replaced by
$${\cal R}^{\half}_N:= e^{-  N^{-1}(P + i  Q)_{\half}} {\cal R}_N e^{  N^{-1} ( P + i  Q)_{\half}} \in \Psi^2_N(\R^1\times
\R).$$
We expand in powers of $N$ to get:
\begin{equation} {\cal R}^{\half}_N\sim \sum_{n=o}^{\infty} N^{-4 + n} \sum_{j+m=n} \frac{i^j}{j!}
(ad(P + i Q)_{\half} )^j {\cal R}_{2 - \frac{m}{2}}\end{equation}
$$:= N^{-4} + N^{-2} {\cal R} +\sum_{n=3}^{\infty}
 N^{-4 + n}{\cal R}^{\half}_{2 - \frac{n}{2}}.$$ 
An obvious induction as in \cite{Z.1} gives that
$$ad( P + i Q)_{\half})^j {\cal R}_{2 - \frac{m}{2}} \in C^{\infty}(\R, {\cal E}_{\epsilon}^{m+j-4}){\cal R}^2
+ C^{\infty}(\R, {\cal E}_{\epsilon}^{m+j-2}){\cal R} + C^{\infty}(\R, {\cal E}_{\epsilon}^{m+j}).$$
It follows that
${\cal R}^{\half}_{2 - \frac{n}{2}}$ has the same filtered structure  as 
${\cal R}_{2 - \frac{n}{2}}.$

We now conjugate ${\cal R}^{\half}_N$ with an exponential of the form $e^{N^{-2}  (P + i Q)_1}.$    As
above, it leads to  a boundary problem  analogous to the previous  case with $j = 1$ except
that now    diagonal terms  with $m = n$ do occur.  In the elliptic (resp. hyperbolic) case, a  diagonal term is a function of $|z|^2$ (resp. $|w|^2$), hence is even
under the involution $z \to \bz$ (resp. $w \to \bar{w}$).   Since the equations for the even/odd coefficients
decouple completely,  the boundary condition reduces to:  
\begin{equation}\begin{array}{l} P^d_{1}(0, |z|^2) = P^d_{1}(L, |z|^2) = 0\\ \\  Q^d_{1}(0, |z|^2) = Q^d_{1}(L, |z|^2) = 0. \end{array}\end{equation}
Similarly in the hyperbolic case with $w$ in place of $z$. 

The same is true for any even $j$ so let us consider the general elliptic case.
We write the diagonal terms in the form:
\begin{equation}\begin{array}{l} P_{j}^d(\s, |z|^2) = \sum_{m: m \leq j} p_{j m}^d(\s) |z|^{2m} \\ \\ 
Q_{j}^d(\s, |z|^2) = \sum_{m: m \leq j} q_{j m}^d(\s) |z|^{2m} 
\end{array} \end{equation}
It is obviously impossible to solve the boundary problem:
\begin{equation}  \begin{array}{l}
\frac{d}{ds} p_{j m}^d(s)   =   \Im {\cal R}_{2 - j }^d|_o \\ \\
 \frac{d}{ds} q^d_{jm}(s)  = - \Re {\cal R}_{2 - j}^d|_o  \end{array} \end{equation}
with zero boundary conditions in general.  To satisfy the boundary condition we need
to add  terms $f_j(|z|^2)$ to the right side.  Then we get:
\begin{equation} \label{DC} \begin{array}{l} \frac{d}{ds} p^d_{j m}(s)   = \Im  {\cal R}_{2-j}^{j, d}|_o - \Im f_j(|z|^2),\;\;\;\;\; 
\frac{d}{ds} q^d_{j m}(s)   = - \Re  {\cal R}_{2- j}^{j, d}|_o + \Re f_j(|z|^2) \end{array} \end{equation}
with the boundary condition above on $p^d_{j m}, q^d_{j m}.$ 
We can solve the equations with 
\begin{equation}\begin{array}{l} p^{d}_{j m}(s)   = \int_0^{\s} \{ -\Im  {\cal R}_{2 -j}^{j, d}|_o(u, |z|^2) + \Re f_j (|z|^2)\} \\ \\
q^{d}_{j m}(s)   = \int_0^{\s}
 \{  \Re {\cal R}_{2-j}^{j, d}|_o(u, |z|^2) + \Im  f_j (|z|^2)\}  du, \end{array} \end{equation}
where
\begin{equation} f_j(|z|^2) = -  \frac{1}{L} \int_0^{L}  {\cal R}_{2-j}^{j, d}|_o (u, |z|^2) du. \end{equation}
In the case $j = 1$, we have then conjugated ${\cal R}_N^+$ to $N^4 + N^2 {\cal R} 
+ Op^w(f_1(I)) + O(N^{-1}).$ We note that $Op^w(f_1(I))$ is a function of $\hat{I}$, so
we  have conjugated to normal form to fourth order.

We then proceed inductively to define   polynomial symbols
$(P + i Q)_{\frac{j}{2}}(\s, \y, \bareta), $  polynomials $f_j(\hat{I})$ and
unitary  $N$-pseudodifferential operators
$W_{N \frac{j}{2}}:=exp( N^{-j} (P + iQ)_{\frac{j}{2}})$ 
such that:

\noindent$$\begin{array}{l}(i)\;\;{\cal R}_{N}^{ \frac{K}{2}}:=
W_{N \frac{K}{2}}^{-1} W_{N \frac{K-1}{2}}^{-1}
 \dots W_{N \frac{1}{2}}^{-1} {\cal R}_N
W_{N \half} \dots W_{N \frac{K-1}{2}}W_{N \frac{K}{2}};\\

(ii)\;\;{\cal R}_N^{\frac{K}{2}} \sim N^{4} + N^{2} {\cal R} +
\sum_{n=3}^{2K} N^{4 - n}{\cal R}^{\frac{K}{2}}_{2 - \frac{n}{2}} + E_{N \frac{K}{2}}^{+ } ;\\

(iii)\;\mbox{ for}\; K \geq n-2,
{\cal R}^{\frac{K}{2}}_{2 - \frac{n}{2}} = {\cal R}^{\frac{n-2}{2}}_{2 -\frac{n}{2}}
=[{\cal R}, (P + iQ)_{\frac{n}{2}-1}] + {\cal R}^{\frac{n-3}{2}}_{-\frac{n}{2}}\\

(iv)  \;\; {\cal R}^{\infty}_{2-j}(\s, D_{\s}, \y, D_{\y})|_o = {\cal
R}^{\infty,o}_{2-j}(\s, \y, D_{\y})|_o =
  f_j(\hat{I})|_o \\

(v) \;\; {\cal R}^{\infty}_{2-\frac{2K+1}{2}}
(\s, D_{\s}, \y, D_{\y})|_o = {\cal R}^{\infty, o}_{2-\frac{2K+1}{2}}(\s,\y,
D_{\y})|_o = 0;\\

 (vi)\;\;{\cal R}^{\frac{K}{2}}_{2-\frac{m}{2}} \in C^{\infty}(\R,
{\cal E}^{m-4}_{\epsilon}){\cal R}^2
 +C^{\infty}(\R, {\cal E}^{m-2}_{\epsilon}){\cal R} +C^{\infty}(\R, {\cal E}^{m}_{\epsilon});\\

(vii)\;\; (P + iQ)_{\frac{K}{2}} \in C^{\infty}(\R,{\cal E}^{\frac{K}{2}+1}_{\epsilon});\\

(viii)\;\;(P^o + iQ^e)_{\frac{K}{2}}(0) = (P^o + iQ^e)_{\frac{K}{2}}(L) = 0; \\

(ix) \;\;E_{N \frac{K}{2}}^+\;\mbox{is}\;\mbox{divisible}\;\mbox{by}\;N^{4-(2K+1)}.
\end{array}$$

 The details about the degrees and parities of the polynomials are  similar to the boundaryless case of \cite{Z.1}
and we therefore omit the details.  An important point which we must address here is that
the real parts $P_{j/2}$ do not contribute to the principal part of the exponent
$(P + i Q)_{j/2}$.  More precisely, the homogeneous part of $P_{j/2}$ of leading order
$2j+2$ equals zero.  This happens essentially because the Weyl symbols of $\bDelta$,
$\bDelta_N$  and ${\cal R}_N^+$
are real,   all  three operators being self-adjoint  with respect to Lebesgue measure.  Conjugation by
unitary operators $e^{i N^{-j} Q_{j/2}}$  preserves this reality. At first sight it appears  that the ${\cal R}_{2 - \frac{m}{2}}^{\frac{j}{2}}$ operators continue to have real symbols
and therefore that $P_{j/2}$ is zero.   However, the $|_o$ operation kills part of
the operator and in particular it kills the self-adjointness and reality of symbols. The ${\cal R}_{2 - \frac{m}{2}}^{\frac{j}{2}}$ generally  have complex
symbols and one does needs to conjugate with $e^{ N^{-j}(P_{j/2} + i Q_{j/2})}.$

What is true is that the leading order terms of order $m$ the Weyl symbol of  ${\cal R}_{2 - \frac{m}{2}}^{\frac{k}{2}}$ in $(y, \eta)$ is real.  As is visible from (vi) aboves, such terms have no factor of ${\cal R}$ in front
while any term with a factor of ${\cal R}$ has a lower degree in $(y, \eta).$  Hence  taking $|_o$ commutes with setting
$\sigma_{{\cal R}} = 0$ in terms of top degree. But setting $\sigma_{{\cal R}} = 0$ in a real symbol produces another
real symbol, so the terms of top degree in $(y, \eta)$ must be real.

Let us  consider the statement (ix) about the error terms.  In each conjugation we  expand the exponential up
to order $N^{-K}$ 
and  absorb the remainder  in the higher powers of $N^{-1}.$  

 To analyse the error, we write
the partial Taylor expansion with remainder of the exponential function as:
$$\begin{array}{ll} e^{ix} = e_M(ix) + r_M(ix),  & e_M(ix) = 1 + ix +
\dots+ \frac{(ix)^M}{M!} \\ & \\
r_M(ix)=(ix)^{M+1}b_M(ix),  &  b_M(ix) = \int_o^1 \dots \int_o^1 t_M^M t_{M-1}^{M-1}
 \cdots t_o^o e^{t_M t_{M-1} \cdots t_o ix}
dt_o \cdots dt_M. \end{array}$$
Clearly $b_M$ is a bounded function on the axis of real $x$.  Now plug in $i x =  N^{-j} {\bf ad}(P + i Q)_{j/2})$; although it is complex, its principal part is real. We get 
that \begin{equation}\begin{array}{l} e^{  N^{-j} {\bf ad}(P + i Q)_{j/2}} = I +  N^{-j} {\bf ad}(P + i Q)_{j/2} +
\dots+ \frac{(N^{-j} {\bf ad}(P + iQ)_{j/2}))^M}{M!}\\ \\ + 
[N^{-j}  {\bf ad}(P + iQ)_{j/2}] ^{M+1} b_M( N^{-j} {\bf ad}(P + iQ)_{j/2}). \end{array} \end{equation}

At the $2K$th stage we are conjugating with  
$2K$ factors $e^{i N^{-j} Q_{j/2}}, j = 1, \dots, K.$
To get a normal form up to order $N^{4 - (2K + 1)}$ at the Kth stage it suffices to choose each $M_j$ such that $j (M_j + 1) > K.$  The $e_M$ terms give  a polynomial differential
operator satisfying (vi) and the $b_M$ terms give the  error $E_{N \frac{K}{2}}^{+}$.  By construction, it is  $N^{4-(2k + 1)}$ times  a  conjugate of a polynomial differential operator of degree $K +2$ in $(\y, \eta).$   \qed

\subsubsection{Remarks}

\noindent{\bf (a)} The operators $\tilde{W}_k^+$ with which we began are now determined by (\ref{W}), i.e.
by $\tilde{W}_k^+ = T_{N_k} W_k^+ T_{N_k}^{-1}.$   We have:
\begin{equation} \tilde{W}^+_k = \Pi_{j = 1}^{\infty} e^{N^{-j} T_{N_k} (P + i Q)_{j/2} T_{N_k}^{-1}}. \end{equation}
As will be clear below the $\tilde{W}_k$'s belong to the class of homogeneous Fourier integral operators, unlike
the $W_k^+$'s, which are exponentials of isotropic pseudodifferential operators.

\noindent{\bf (b)} In the next section we will glue together the $\tilde{W}_k^+$'s's
into a Fourier integral operator.  The vanishing of the principal part of $P_{j/2}$ then implies
that it does not contribute to the principal symbol of $(P + i Q)_{j/2}$, but only to the
amplitude. It will also imply that $b_M( N^{-j} {\bf ad}(P + iQ)_{j/2})$ is a bounded FIO.

\noindent{\bf (c)} Finally, let us note that $C W_N^+ e^{i N^2 L} C = W_N^+ e^{i N^2 L} $ if and only if $e^{- i N^2 L } = e^{i N^2 L}$, i.e.
if and only if $ N^2  = \frac{\pi}{L} k$ for some $k$.  Henceforth we put $N_k = \sqrt{\frac{\pi}{L} k}$ and
only consider these special values of $N$. For notationally simplicity we write $W_k^+$ for
$W_{N_k}^+$.  We also define  $F(k, \hat{I})$ as the formal asymptotic series whose
$K$th partial sum is  $ F_{N_k K}(\hat{I})$.

\section{Conjugation of the Dirichlet Wave group}

So far, we have conjugated the semiclassical Laplacian on the open space in a microlocal
neighborhood of $\gamma$ to a semiclassical normal form modulo a small semiclassical 
remainder. We now glue together the component intertwining operators $\tilde{W}_k^+ = T_{N_k} W_k^+ T_{N_k}^{-1}$  into a homogeneous Fourier integral intertwining operator of the Dirichlet wave group to its normal form. We emphasize that the intertwining operators are defined in
the open space.

Conjugation by the scaling operator $T$ and its local form $T_{N_k}$ plays an important role since it converts 
isotropic pseudodifferential operators into standard homogeneous pseudodifferential operators.  This will be
discussed further below.  Since it becomes heavy notationally to distinguish the action operators $\hat{I}$ from
their conjugates we abuse notation somewhat by denoting both $\hat{I}$ and $T \hat{I} T^{-1}$ by the same
symbol $\hat{I}$. 

By the Dirichlet wave kernel we mean the fundamental solution $E(t, x, y)$ of the mixed
wave equation (\ref{DWK}) with Dirichlet boundary conditions, i.e.  the kernel
of $\cos t \sqrt{\bDelta_{\Omega}})$.  The `free' Laplacian $\bDelta$ on the open space
is only well-defined in a neighborhood of $\gamma$, so to be precise we need to cut it off
to $U_{\epsilon}$.  We will not indicate this cutoff in the notation because we will
explicitly microlocalize it later on.

\subsection{Local components on $S^1_{2L} \times \R$}

By `local component' we mean an operator which only acts on a specific Fourier coefficient.

\begin{defn} Define the operators:   $\tilde{W}_k: L^2(S^1_{2L} \times \R) \to L^2(S^1_{2L} \times \R)$

$$\begin{array}{l} \tilde{W}_k  = \tau_{\epsilon} (\tilde{W}_k^+ \Pi_{+k} - \tilde{W}_k^- \Pi_{-k}).\end{array}$$
Here, $\tilde{W}_k^- = C \tilde{W}_k^+\Pi_{+k}  C.$ 
\end{defn}

In 

\begin{prop}\label{LOC} We have: 

\noindent(i) 
$$ \tau_{\epsilon} \bDelta   \tilde{W}_k^{\pm} \Pi_{\pm k} = \tau_{\epsilon}  \{\tilde{W}_k^{\pm} F(|D|, \hat{I})^2\Pi_{\pm k} + \tilde{W}_k^{\pm} E_{k K}^{\pm } \Pi_{\pm k}\}.$$
Here, $|D| = \sqrt{- D^*D}$ on $S^1_{2L}.$

\noindent(ii)
$$ \tau_{\epsilon} \bDelta   \tilde{W}_k \Pi_o = \tau_{\epsilon}  \{\tilde{W}_k F(|D|, \hat{I})^2\Pi_o + \tilde{W}_k E_{k K}^+\Pi_o\}.$$

\noindent(iii) The Schwarz kernel of $\tilde{W}_k \Pi_o$ vanishes on $\partial \Omega_o.$
\end{prop} 

\noindent{\bf Proof} 

\noindent(i)  This is essentially just a restatement of Lemma (\ref{MAIN}). 
We can prove it by testing both sides against functions of the form $e^{i k s} f(N_k y)$. 
  We do the case $k > 0$. On the set  $\tau_{\epsilon} \equiv 1$ we have:
$$\begin{array}{l}  \bDelta  \tilde{W}_k^{+} e^{i \frac{\pi \s k}{L}} f(N_k  y) 
= e^{i \frac{\pi \s k}{L}}  e^{- i \frac{\pi \s k}{L}}\bDelta e^{i \frac{\pi \s k}{L}}  T_{N_k} W_k^+ T_{N_k}^{-1} f(N_ky)\\ \\ 
= e^{i \frac{\pi \s k}{L}}T_{N_k}  \bDelta_{N_k}   W_k^+ f(y)
=  e^{i \frac{\pi \s k}{L}} T_{N_k} W_k^+ (F(k, \hat{I})^2 + E_k^K) f(y) \\  \\
 =  \tilde{W}_k^+ (F(k, \hat{I}))^2  e^{i \frac{\pi \s k}{L}} f(N_k y)  +   \tilde{W}_k^+ E_k^K e^{i \frac{\pi \s k}{L}} f(N_k y) \\ \\
=  \{\tilde{W}_k^+ (F(|D|, \hat{I}))^2  +  \tilde{W}_k^+ E_k^K)\} \Pi_{\pm k }e^{i \frac{\pi \s k}{L}} f(N_k y). 
\end{array}$$
Between lines two and three, the action operator $\hat{I}$ changed from its isotropic form to its homogeneous form.
The case of $k <0$ follows by taking complex conjugates.

\noindent(ii)  This follows immediately from (i) by 
taking the difference of the $+$ and $-$ cases.

\noindent(iii) The Schwarz kernel of $\tilde{W} \Pi_o$ is given by $\sum_{k} \tilde{W}_k \Pi_{ok}$, so it suffices
to show that $\tilde{W}  \sin (\frac{k \pi \s}{L}) f(N_k y) = 0$  on $\s = 0, L$. This follows from 
Lemma (\ref{MAIN}), since $(\tilde{W}_k  \sin (\frac{k \pi \s}{L}) f(N_k y) = 0$  on $\s = 0, L$.  Also, we obviously
have $ \sin (\frac{k \pi \s'}{L}) f(N_k y') = 0$  on $\s = 0, L$. 
\qed

\subsection{The homogeneous intertwiner}

We glue the $\tilde{W}_k$ together as follows:

\begin{defn} 
Define the operators $\tilde{W}^{\pm}, \tilde{W} : L^2(S^1_{2L} \times \R_{\y}) \to L^2(S^1_{2L} \times \R_{\y})$ by
$$\begin{array}{ll} (i) & \tilde{W}^+ = \sum_{k = 1}^{\infty} \tilde{W}_k^+  \Pi_k^+  \\ & \\
(ii) & \tilde{W}^- =   \sum_{k = 1}^{\infty}    \tilde{W}_k^- \Pi_k^-  \\ & \\
 (iii) & \tilde{W} = \tilde{W}^+ - \tilde{W}^-.
\end{array} $$
\end{defn}

 We now analyse the Fourier integral nature of $\tilde{W}, \tilde{W}^{\pm}$. A key point is that  the operators $T_{N_k}$ glue together to form a `global' transverse scaling operator
\begin{equation}\begin{array}{l} 
T : L^2([0, L]_{\s} \times \R_{\y}) \to L^2([0, L]_{\s} \times \R_{\y}) \\ \\  
  T \sin ( \frac{\pi k \s}{L}) f(\y) : = N_k \sin (\frac{\pi k \s}{L}) f(N_k y) .
\end{array} \end{equation}
Conjugation with $T$ transforms the isotropic Weyl calculus into the usual homogeneous psuedodifferential calculus.
   For instance 
\begin{equation} T (D_{\y}^2 + \y^2) T^* = |D_s|^{-1} D_{\y}^2 + \y^2 |D_s|,\;\;\;\;\;
T (D_{\y}^2 - \y^2) T^* = |D_s|^{-1} D_{\y}^2 - \y^2 |D_s|. \end{equation}
As is verified in \cite{G} and elsewhere, $T$ is an oscillatory integral operator associated
to the canonical transformation
\begin{equation} \psi(s, \sigma, y, \eta) = (s + \frac{y \eta}{2 \sigma}, \sigma,
\sqrt{\sigma} y, \sqrt{\sigma}^{-1} \eta). \end{equation}

The following proposition is  analogous way to 
(\cite{Z.1}, Proposition (3.4)). Below, the notation $A \sim B$ in $V$ means that $A - B$ is smoothing
in $V$.  In the following proposition we assume the order $K$ of the Birkhoff normal form
is infinity. 

\begin{prop}\label{FIO} There exist conic neighborhoods $V^{\pm}$, $V^{\pm '}$ of $\gamma_o^{\pm}:= \gamma_o \cap \{ \pm \sigma > 0\}$   and 
 canonical transformations  $\chi^{\pm} : V^{\pm} \to V^{\pm '}$ such that:

\noindent(i)  $\chi^{\pm} = Id$ on   $\R \gamma_o^{\pm} := \; \{\s, \sigma, \y,
\bareta) \} \in V^{\pm}: \y = \bareta = 0\}$;

\noindent(ii)  $  \tilde{W}^{\pm}  \in  I^0(S^1_{2L} \times \R \times S^1_{2L} \times
\R, \mbox{gr}\chi^{\pm})$
where $\mbox{gr} \chi^{\pm} \subset V^{\pm} \times V^{\pm '}$ is the graph of $\chi^{\pm}$;

\noindent(iii) $ \tilde{W}^{\pm} $  is elliptic in $V^{\pm} \times V^{\pm  '}$.
\end{prop}

\noindent{\bf Proof} 

To analyse the  sum over $k$, we
 enlarge the  `open model space' $S^1_{2L} \times \R$ to the `space-time' $S^1_{2L} \times \R \times S^1_{2L}$
and define   the operator 
\begin{equation}\begin{array}{lll}  e^{i t |D_{\s}|}& : 
  L^2(S^1_{2L} \times \R) \to L^2(S^1_{2L} \times \R \times S^1_{2L}), & 
e^{i t |D_{\s}|} e^{i \frac{i \pi k \s}{L}} f(\y)   =
e^{i \frac{i \pi k t}{L}} e^{i \frac{i \pi k \s}{L}} f(\y). \end{array} \end{equation}
The range of $e^{i t |D_{\s}|}$ is contained in the kernel ${\cal H}$ of the `wave
operator' $ |D_{\s}|^2 - |D_t|^2$ on $S^1_{2L} \times \R \times S^1_{2L}$. We
also denote by $P$ the orthogonal projection to ${\cal H}$. 
 We then introduce the further operators:

\begin{defn} Define $ \tilde{W}_{|D_t|}^{\pm} (\s, \y, D_{\y})$ and $\tilde{W}_{|D_t|}$ on $L^2(S^1_{2L} \times \R \times S^1_{2L})$

\noindent(a)  For $k \geq 0$ put 
$$ \tilde{W}_{|D_t|}^+(\s, \y, D_{\y}) \Pi^+ f(\s, N_k \y) e^{i \frac{\pi k t}{L}} = e^{i \frac{\pi k t}{L}} \tilde{W}^+_k  f(\s, N_k \y);$$
For $k < 0$ put it equal to zero. 

\noindent(b) $\tilde{W}_{|D_t|}^- (\s, \y, D_{\y}) = C \tilde{W}_{|D_t|}^+(\s, \y, D_{\y}) C; $

\noindent(c)  $\tilde{W}_{|D_t|} (\s, \y, D_{\y}) =  \tilde{W}_{|D_t|}^+ (\s, \y, D_{\y}) \Pi_+ -  W_{|D_t|}^- (\s, \y, D_{\y})\Pi_-.$
\end{defn}

These definitions suggest the introduction of a scaling operator $\tilde{T}$ adapted to the $t$-variable.
We therefore define  $\tilde{T} e^{it m} f(s,  y)  = N_m e^{i t m} f(s, N_m y).$ We are only interested in
its action on the invariant subspace 
${\cal H}$ where the frequency in $s$ and $t$ are the same.  As with $T$, $\tilde{T}$ is an oscillatory integral operator with underlying canonical transformation
$$\tilde{\psi}(s, \sigma, t, \tau, y, \eta) = (s, \sigma, t + \frac{y \eta}{2 \tau},
\tau, \sqrt{\tau} y, \sqrt{\tau}^{-1} \eta).$$
We then have: 
$$\tilde{W}_{|D_t|}^+(\s, y, D_y)   = \Pi_{j = 0}^{\infty} e^{|D_t|^{-j/2} (\tilde{P}+ i  \tilde{Q})_{j/2}(\s, \y, D_{ \y})}, \;\;\;\;\mbox{with}\;\;\tilde{P} = \tilde{T} P \tilde{T}^{-1}, \; \tilde{Q} =  \tilde{T} Q \tilde{T}^{-1}. $$
The operators $\tilde{P}, \tilde{Q}$ are (usual) homogeneous pseudodifferential operators, as we will argue below.

The following identity is the key to the Fourier integral properties of $W^{\pm}, W$: Let
$j: S^1_{2L} \times \R \to S^1_{2L} \times S^1_{2L} \times \R $ be the inclusion
$j(s, y) = (s, 0, y).$  Then:

\begin{equation} \begin{array}{l}\tilde{W}^{\pm} (\s, \y, D_{\y}) = j^* \tilde{W}_{|D_t|}^{\pm} (\s, \y, D_{\y})   e^{i t |D_{\s}|}\\ \\
\tilde{W} = j^* \tilde{W}_{|D_t|}   e^{i t |D_{\s}|}\end{array} \end{equation}
To prove it,   we apply both sides to functions of the form  $e^{i \frac{\pi k \s}{L}} f(N_k y)$. In the plus case,
we have (for $k \geq 0$),
$$\begin{array}{l} j^* \tilde{W}_{|D_t|}^+(\s, \y, D_{\y})   e^{i t |D_{\s}|} e^{i \frac{\pi k \s}{L}} f(N_k y)\\ \\ = j^* e^{i \frac{\pi |k| t}{L}} \tilde{W}_{k}^+(\s, \y, D_{\y}) e^{i \frac{\pi k \s}{L}} f(N_k y) \\ \\
= \tilde{W}^+ e^{i \frac{\pi k \s}{L}} f(N_k y) \end{array}$$
and similarly in the minus case.  

 We now complete the proof of (i)--(iv). For simplicity we only consider the $+$ case,
the other being essentially the same. By definition,    $\tilde{W}_{|D_t|}^+$ is a product of factors of the form $e^{ |D_t|^{-j/2}(\tilde{P}_{j/2} + i \tilde{Q}_{j/2})}$ with $\tilde{P}_{j/2}, \tilde{Q}_{j/2}$ 
scaled polynomial differential operators of degree $j + 2$  in
the variables $|D_t|^{\half} y$ and $|D_t|^{-\half} D_y$. As noted in the proof of Lemma (\ref{MAIN}),  the terms in $\tilde{P}_{j/2}$ have the parity of $j$ and the leading order term in $(y, D_y)$ vanishes.  Hence,  $|D_t|^{-j/2}\tilde{P}_{j/2}$ is actually a pseudodifferential  operator of at most zero order and its  
does not affect the status of $\tilde{W}_{|D_t|}^+$ as a Fourier integral operator.
 Since  $\tilde{Q}_{j/2}$ is a polynomial  order $j + 2$ in
($|D_t|^{\half} y, |D_t|^{-\half} D_y$), it follows that  $|D_t|^{-j/2} (\tilde{P} + i \tilde{Q})_{j/2}$ is a first order (homogeneous) pseudodifferential
operator  of real principal type for each $j$.  Its exponential   $ e^{i |D_t|^{-j/2} (\tilde{P} + i \tilde{Q})_{j/2}}$ is therefore  a Fourier integral operator.  Hence the product of  any finite number of factors is a Fourier integral operator.    As discussed above, we only  need
to deal with finitely many factors, but we note that the the full infinite product can be regularized as a Fourier
integral operator.  Indeed,  one  may apply the Baker-Campbell-Haussdorf formula to any finite number of factors to
rewrite  $\tilde{W}^+_{|D_t|}$ in the form  $\tilde{W}^+_{|D_t|} = e^{ \sum_{j = 1}^{\infty} |D_t|^{-j/2}({\cal P}_{j/2} + i {\cal Q}_{j/2})}$. If one interprets the sum as given by  the Borel summation method then the exponent is well-defined
as a pseudodifferential operator of order one and of real principal type. 
 The principal symbol of the series $\sum_j |D_t|^{-j/2} ({\cal P} + i {\cal Q})_{j/2}$ is then well-defined as a 
an element of $S^{1, 2}$ of the form $H(s, y, \eta, \tau):= \sum_j \tau^{-j/2} \tilde{q}_{j/2}(s, \sqrt{\tau} y, \sqrt{\tau}^{-1} \eta)$ in $(y, \eta)$, with $\tilde{q}_{j/2}$ the leading order homogeneous part of the Weyl 
symbol of $\tilde{Q}_{j/2}$.    The symbol of ${\cal P}_{j/2}$ is of order at most zero as a symbol and hence does
not contribute to $H$.

The operators $e^{i t |D_s|}$ and $j^*$ are clearly Fourier integral operators, with
associated canonical relations
\begin{equation}\begin{array}{l}  C = \{ ((s, \sigma,
y, \eta);(s + t, \sigma, t, \sigma, y, \eta))\} \subset T^*(S^1_{2L} \times \R \times S^1_{2L} \times S^1_{2L} \times \R)  \\ \\
\Gamma_{j^*} = \{((s, \sigma, y, \eta); (s, \sigma, 0, \sigma, y, \eta)\} \subset 
T^*(S^1_{2L} \times \R \times S^1_{2L} \times \R \times S^1_{2L} ). \end{array} \end{equation}

The canonical relation underlying $ j^* \tilde{W}^+_{|D_t|} e^{i t |D_s|} $ is therefore
given by the composite relation
$$ \Lambda = \Gamma_{j^*} \circ \; \mbox{gr}\; \phi^1 \circ C $$
where $\phi^u = (\exp u \Xi_H)$ with  $\Xi_H$  the Hamiton vector field of $H$ and $\exp u \Xi_H$  its flow.  
The first two factors compose to the relation $\{((s, \sigma, y, \eta); \phi^1(s + t, \sigma,
t, \sigma, y, \eta))\}.$ It is easy to see that  $\phi^1(s + t, \sigma, t, \sigma, y, \eta)$
has the form $(s + t, \sigma(1), t(1), \sigma, y(1), \eta(1))$ where $x(1)$ stands for
the value at $u = 1$ of the of $x$-coordinate of the orbit of $\phi^u$ through 
the initial point $(s + t, \sigma, t, \sigma, y, \eta).$ To contribute to $\Lambda$
one must have $0 = t(1) = t + \int_0^1 \frac{\partial H}{\partial \tau}(\phi^u (s, \sigma,
t, \sigma, y, \eta) du.$ We observe that $H$ vanishes to order at least two along
 $\R^+ \gamma_o = \{y = \eta = 0\}$ and hence that $\phi^u$ acts as the identity on this set. In particular, $t(1) = t = 0$
there. Futhermore, for $(y, \eta)$ sufficiently small, the derivative of $t(1)$ with
respect to $t$ is non-zero and hence there exists a unique solution $t(s, \sigma, y, \eta)$
of the equation relating $t(1)$ and $t$. It follows that (at least) in a sufficiently small cone
around $\{y = \eta = 0 \}$, $\Lambda$ is the graph
 of the canonical transformation $\chi^+(s, \sigma, y, \eta) = \phi^1(s + t(s, \sigma, y, \eta),
\sigma, t(s, \sigma, y, \eta), \sigma, y, \eta)$ (where the $(t, \tau)$-coordinates are
omitted on the left side).

The $-$ component has a similar description. Therefore, 
 $W$ is a  microlocal Fourier integral operator 
associated to the  graph of a canonical transformation defined in a small cone around
$\R^+ \gamma_o$. 
 Each of the three operators in its composition  has a canonical  principal symbol and
it follows in a standard way that  $\sigma(W)$ is a graph 1/2-density.  Hence $W$ is
microlocally elliptic, concluding the proof.\qed

Below we will need a somewhat smaller open cone with the following property:
\begin{prop}\label{VZERO} Let $int(\gamma_o) = \gamma_o \cap T^*(int \Omega_o)$ where $int \Omega_o$ denotes the interior
of $\Omega_o.$ Then there exists an open conic neighborhood $V_o$ of $ int(\gamma_o)$ with the property that
$\chi(V_o)$ lies in the interior of $T^*(\Omega_o)$. \end{prop}

\noindent{\bf Proof} The point we must address is that $\chi$ could carry a covector at $(s, y) \in int(\Omega_o)$
to a covector on $\partial \Omega_o$ or in the exterior.  To determine whether this happens we must study
$s' = \frac{\partial}{\partial \sigma} H.$ We observe that this component is homogeneous of degree 0 on $T^*\Omega_o)$
and that it vanishes to order three along $y = \eta = 0.$  It follows that there is a neighborhood $V_o$ of
$\gamma_o$ of the form $\max \{|s|, |L - s|\} \leq C \min \{ |y|^3, |y|^2 \frac{|\eta|}{\sigma}, |y| \frac{|\eta|^2}{\sigma^2},  \frac{|\eta|^3}{\sigma^3} \}$ with the property that $\phi(s, y, \sigma, \eta) \in 
T^*(int \Omega_o)$ if $(s, y, \sigma, \eta)  \in V_o$. \qed

\subsection{The error term}

We now make a similar analysis of the error term.
We first define operators $\tilde{E}^{\pm K}$  on  $H^2_{\pm} (S^1_{2L} \times \R) $ and $\tilde{E}_K$ on 
$L^2 (S^1_{2L} \times \R) $
by 
$$ \begin{array}{l} \tilde{E}_K^{+} e^{i k \frac{ s \pi}{L}}f(N_k y)  =   \tilde{W}_k^+ \tilde{E}^{+}_{N_k K} e^{i k \frac{ s \pi}{L}}f(N_k y),\;\;\;\;(k \geq 0) \\ \\
\tilde{E}^{-K} e^{i k \frac{ s \pi}{L}}f(N_k y) = 
 C  \tilde{W}_K^+  \tilde{E}^{+}_{N_k K} C   e^{i k \frac{ s \pi}{L}}f(N_k y),\;\;\;(k < 0) \\ \\
\tilde{E}_K = \tilde{E}_K^{+} \Pi_+ - E_K^{-} \Pi_-\end{array}$$

To analyse the kernels of $\tilde{E}_K^{\pm }$ and $\tilde{E}_K$ we introduce the operator on $S^1_{2L} \times \R \times S^1_{2L}$:
\begin{equation}\begin{array}{l} \label{Eplus}
 \tilde{E}_{|D_t| K}^{ + }  =  \tilde{W}_{|D_t|}^{+ *}  \bDelta \tilde{W}_{|D_t|}^{+ } - 
F_K^{ 2}(|D|, \hat{I}).\end{array} \end{equation}
Recall here that $F_K^{ 2}(|D|, \hat{I}) = \sum_{j = 1}^{2K} f_j(|D_t|, \hat{I})$.

\begin{prop}  As an  operator on $L^2(S^1_{2L} \times \R)$, we have: $\tilde{E}_K^{+} = j^* \tilde{W}_{|D_t|}^+  \tilde{E}_{|D_t|K}^{ + }   e^{i \frac{\pi |D_{\s}| t}{L}}.$ \\

\end{prop}

\noindent{\bf Proof}\\

  As above, we test both sides on functions of the form $e^{i k \frac{ s \pi}{L}}f(N_k y)$ for $k>0$:
$$\begin{array}{l} j^* \tilde{W}_{|D_t|}^+ \tilde{E}_{|D_t|}^{+K}  e^{i \frac{\pi |D_{\s}| t}{L}}   e^{i k \frac{ s \pi}{L}}f(N_k y)
= j^* \tilde{W}_{|D_t|}^+ \tilde{E}_{|D_t|}^K  e^{i \frac{\pi k t}{L}} e^{i k \frac{ s \pi}{L}}f(N_ky)\\ \\
= 
 j^* W_k^+ E_{N_k K}^{+ } e^{i \frac{\pi k t}{L}} e^{i k \frac{ s \pi}{L}}f(N_k y)\\ \\
= 
  W_k^+ E_{N_k K}^{+ }  e^{i k \frac{ s \pi}{L}}f(N_k y) = \tilde{E}_K^{+ }  e^{i k \frac{ s \pi}{L}}f(N_k y). \end{array}$$ 
Similarly for $k < 0$. \\
\qed

In the following we denote by $\tilde{V}$ the cone in $T^*(S^1_{2L} \times \R \times S^1_{2L})$
defined by $(s, \sigma, y, \eta, t, \tau) \in \tilde{V}$ iff $(s, \sigma, y, \eta) \in V,
C < |\tau|/ \sigma \leq 1/C$ for some $C < 1.$  We also denote by $\R^+ \gamma \times
T^*S^1_{2L}$ the symplectic symplectic subcone of $\tilde{V}$ which is defined by
$y = \eta = 0.$

\begin{prop} \label{ERROR} $ \tilde{E}_{|D_t| K}^{+}  \in Op S^{2, K}(\tilde{V}, \R^+ \gamma \times T^*S^1_{2L}).$\end{prop} 

\noindent{\bf Proof}  By (\ref{Eplus}),  by Proposition (5.6) and by Egorov's theorem, we see that $\tilde{E}_{|D_t| K}^{+}$ is  a homogeneous
pseudodifferential operator of order $2$.   To show that its symbol lies in
$S^{2, K}$ it suffices  by proposition (\ref{SYM}) to show
that the Taylor expansion of its scaled symbol is divisible by $N^{-(K + 1)}.$ The
symbol of $ \tilde{E}_{|D_t| K}^{+}$ equals $E^{+ }_{N K}$ with $|\tau|$  substituted
for $N^2$ and we will denote it by $E_{|\tau| K}^{+ }(s, \y, \eta).$  Here we use that
$E^{ +}_{N K}$ is independent of $\sigma$.  Rescaling in $\tilde{V}$ replaces $|\tau|$ by
$N^2$.  Moreover, the transversal scaling $(\frac{1}{N} y, N \eta)$ is canceled by
the operator $T$. Hence the scaled symbol is precisely the original semiclassical symbol 
 $ E_{N K}(s,   \y, \bareta)$.    By lemma (\ref{MAIN}) it is divisible by $N^{-(K + 1)}.$  \qed 

From the previous two propositions we have:

\begin{cor}  ${\bf E}_K$ is a sum of terms of the form $A R B$ where $A, B$
are bounded Fourier integral operators and where $R \in Op S^{2, K}(V, \R^+ \gamma).$
\end{cor}

\subsection{Conclusion}

 The following proposition sums up the discussion.  It gives the homogeneous analogue
of Proposition (\ref{LOC}) and the proof is essentially the same.  We use the notation
$A \sim B$ in $V$ to mean that $A - B$ is a Fourier integral operator  of order $- \infty$ in $V$.
Also we emphasize that the action variables are the scaled ones, i.e. functions of $|D|^{\half} y$ and
$|D|^{-\half} D_y.$

\begin{prop}\label{GLOB} We have:\\

\noindent(i) $  \bDelta  \tilde{W} \sim    \tilde{W} [ F(|D|, \hat{I})^2 +  {\bf E}_K)] $ in $V$.

\noindent(ii) If $u \in H^1_0(\Omega_o)$ then $\tilde{W} u = 0$ on $\partial \Omega_o.$

\noindent(iii) The canonical transformation $\chi$
conjugates the Hamiltonian $|\xi|_{g_{\Phi}}$ on $V^{\pm}$ to the Hamiltonian $F(\sigma, I)$ 
on $V_o^{\pm}$ modulo an error which vanishes to infinite  order on $\R^+ \gamma_o$.

\end{prop}

\noindent{\bf Proof} 

\noindent(i) Since $\tilde{W} = \sum_k (\tilde{W}_k^+  \Pi_{+k} - \tilde{W}_k^-  \Pi_{-k})$ it  follows as in Proposition (\ref{LOC}) that 
$ \tau_{\epsilon} \bDelta  \tilde{W} e^{i k \frac{ s \pi}{L}}f(N_k y) = \tau_{\epsilon} \tilde{W} [ F(|D|, \hat{I})^2 + \tilde{E}_K] e^{i k \frac{ s \pi}{L}}f(N_k y).$

\noindent(ii) This follows by taking the principal symbol of the equation in (i).  In
the principal symbol, the remainder is homogeneous of order 1 and vanishes to order $K$
at $\gamma_o.$

\noindent(iii) It suffices to show that $W \sin k \frac{\pi \s}{L} f(N_k y) = 0$ on $\partial \Omega_o.$
But this follows from the fact that $W_k \sin k \frac{\pi \s}{L} f(N_k y)= 0$ on $\partial \Omega_o$
for all $k$.

\qed

\subsection{Normal form and microlocal parametrix}

We now use  $\tilde{W}$ to conjugate the (odd part of the) wave group of the normal form to a
kind of parametrix for the mixed problem on $\Omega_o.$ 
First, let us define the odd Dirichlet normal form wave group:

\begin{defn} By  $F_o(t, x, y)$ we denote the odd part of the fundamental solution
of the wave equation:
$$\label{NFWK} \begin{array}{ll} \frac{\partial^2 F_o(t)}{\partial t^2} =
F(|D|, \hat{I})^2 F_o(t) & \mbox{on} \; \R \times (S^1_{2L} \times \R) \times (S^1_{2L} \times \R) \\ & 
\\ F_o(0) = \Pi_o & \frac{\partial F_o}{\partial t}(0) = 0 \\
& \\
F_o(t,\cdot,z) = 0 = F_o(t, z, \cdot) & \mbox{on}\;\; \partial \Omega_o. \end{array}$$ \end{defn}

To be correct,  $F$ is only defined in a microlocal neighborhood of $\gamma_o$. We
do not indicate this in the notation since we will microlocalize it later on. 
In the elliptic case, for instance, $F_o$ has the eigenfunction expansion
\begin{equation} F_o(t, \s, \y; \s', \y') = \sum_{k,q } \cos t (F(k, q + \half)) \phi_{kq}(\s, \y) \phi_{kq}(\s', \y')\end{equation}
where $\phi_{kq}(\s, \y) = \sin k \frac{\pi \s}{L} D_q(N_k y)$. 
By odd part we refer to the projection under $\Pi_o$. 

We now wish to conjugate the odd normal form wave group to a microlocal parametrix for the Dirichlet wave group $E(t)$.  

\begin{prop} There  exists a conic neighborhood $V' \subset V$ of $\gamma_o$ and a  zeroth order pseudodifferential operator  $G$ such that:
$$\left \{ \begin{array}{l} \tilde{W} \tilde{W}^* G \sim  \tilde{W}^* G \tilde{W} =  I + R \;\; \mbox{with}\;\; WF(R) \cap V' = \emptyset\\ \\
\tilde{W} F^2 \tilde{W}^*G = \Delta  + A_1 {\bf E}^K A_2 + R_1,\;\;\;\; WF(R_1) \cap V' = \emptyset\end{array} \right.$$
where $A_1, A_2$ are zeroth order Fourier integral operators. \end{prop}

\noindent{\bf Proof}: By proposition (\ref{FIO}), 
 $\tilde{W}$ is a zeroth order Fourier integral operator microsupported in $V \times \chi(V)$ and associated
to the graph of the canonical transformation $\chi$.  Hence   $\tilde{W}  \tilde{W}^*$ is a zeroth order
pseudodifferential operator microsupported in $V$ with symbol identically equal to one
in a smaller cone $V' \subset V$.  Hence
there exists a positive zeroth order self-adjoint  pseudodifferential operator $G $ 
microsupported in $V$ such that $\tilde{W}  \tilde{W}^* G  I+ R$ with $WF(R) \cap V' = \emptyset$. 

Let us put $\tilde{W}_{-1} = \tilde{W} G$.  Then by proposition (\ref{GLOB}), 
\begin{equation} \begin{array}{l} \tilde{W} F^2 \tilde{W}_{-1}    
=   \tilde{W} F^2  \tilde{W}^* G \\ \\
\sim \bDelta \tilde{W}  \tilde{W}^* G + \tilde{W} {\bf E}^K  \tilde{W}^* G \;\; \mbox{in}\;\; V \\ \\
\sim  \bDelta  + R\;\; \mbox{in}\;\;V' \end{array} \end{equation}
where 
$ R =  \tilde{W} {\bf E}^K   \tilde{W}^* G.$
\qed

We now introduce a kind of  microlocal parametrix for the wave kernel.  First we define
a microlocal cutoff $\psi$  to the conic neighborhood $V$ of $\R \gamma_o$, with $\psi \equiv 1$
in a conic neighborhood of $int \R^+ \gamma.$  
\begin{defn} The parametrix is defined by:    ${\cal E}(t) = \tilde{W} F_{oT}(t)\psi_{V} 1_{\Omega_o} \tilde{W}_{-1}$. \end{defn}
The extent to which ${\cal E}(t)$ is indeed a microlocal parametrix near $\gamma_o$ is given in the following
lemma:

\begin{lem}\label{PARA} For any $A \in \Psi^*(\Omega)$ with microsupport  $WF(A) \subset V_o$
(defined in proposition (\ref{VZERO}), 
 $E (t) A $ and ${\cal E}(t)  A $ are Fourier integral operators in the class $I^0((S^1_{2L} \times \R)
\times (S^1_{2L} \times \R), \mbox{gr} \chi) + I^o ((S^1_{2L} \times \R)
\times (S^1_{2L} \times \R),  \mbox{gr} \chi \circ r)$.  We have
$$ \half E (t) A  \sim {\cal E}(t) A   + R(t) A $$
where  $R(t)$ is a finite  sum (or smooth compactly supported  integrals) of terms  $ B_1 R_1(t) B_2$ where $B_1, B_2$ are bounded Fourier integral operators and $R_1\in Op S^{2, K}. $ \end{lem}

\noindent{\bf Proof}: $E (t) $ is uniquely characterized as the microlocal solution of the Cauchy problem (\ref{DWK}) 
on $\Omega_o$ with initial condition equal to the identity operator on $L^2(\Omega_o)$.  Hence 
it suffices to show that ${\cal E}$ is a also a  microlocal solution  modulo errors in the stated class. 

 We first verify that  ${\cal E}t)$
is a microlocal  solution of the wave equation in $V'$.  But
\begin{equation}\label{NFWE} \begin{array}{l} \frac{\partial^2 }{\partial t^2}  
{\cal E} (t)  = \tilde{W} F(|D|, \hat{I})^2 F_{o } \psi_V 1_{\Omega_o} \tilde{W}_{-1}   \\   \\ 
=  \tilde{W} F(|D|, \hat{I})^2 \tilde{W}_{-1} \tilde{W}  F_{o }\psi_V 1_{\Omega_o} \tilde{W}_{-1}    + \tilde{W} F(|D|, \hat{I})^2 
(I -\tilde{W}_{-1} \tilde{W} )  F_{o } \psi_V 1_{\Omega_o} \tilde{W}_{-1}   \\  \\
=  \bDelta {\cal E}(t)  +  E^K {\cal E}(t)  + R_1(t), \end{array} \end{equation}
where $R_1(t) = \tilde{W} F(|D|, \hat{I})^2 
(I -\tilde{W}_{-1} \tilde{W} )  F_{o } \psi_V 1_{\Omega_o} \tilde{W}_{-1} + R_2 {\cal E} (t)$ with $WF(R_2) \cap V' =
\emptyset.$ By shrinking $V$ if necessary, the first term also has order $- \infty$
in $V$.  Hence ${\cal E}(t)$ is a microlocal solution of the forced wave equation
with forcing term of the form $E^K {\cal E}_{\epsilon}(t) + R_1(t).$ 

Regarding the boundary condition, we observe   that 
 $W \Pi_o(x, y) = 0$ if $x \in \partial \Omega$ and 
so ${\cal E}(x, y) = 0$ if $x \in \partial \Omega.$  

The initial condition is more difficult.  We must
show that $\tilde{W} \Pi_o \psi_V 1_{\Omega_o} \tilde{W}_{-1} A \sim \half A$ for a pseudodifferential operator $A$ microsupported in $V_o$  To this
end we first recall  that $\Pi_o = \half (I - r)$ where $r$ is the reflection $(\s, \y) \to (-\s, \y) $ (with $\s$ mod $2 L$.) Hence
$\tilde{W} \Pi_o  \psi_V 1_{\Omega_o} \tilde{W}_{-1} A = \half \tilde{W} \;\psi_V\;1_{\Omega_o} \;\tilde{W}_{-1} A - \half \tilde{W}\; r \;\psi_V \; 1_{\Omega_o} \tilde{W}_{-1} A$.  It suffices
to show that $ \tilde{W} \psi_V  1_{\Omega_o} \tilde{W}_{-1} A \sim A$ and that $ \tilde{W} \; r \; \psi_V \; 1_{\Omega_o} \tilde{W}_{-1} A \sim 0.$ 

The operator  $ \tilde{W} \psi_V 1_{\Omega_o} \tilde{W}_{-1}$ is a pseudodifferential operator with a singular symbol.  Temporarily
ignoring the singularity along $\partial \Omega_o$ we may apply Egorov's theorem formally to find that the complete
symbol is identically equal to one in $V_o$. To see this, it is convenient (although not necessary) to use that $\tilde{W}= e^{i P}$ for some first order pseudodifferential $P$ of real principal type (see proposition (\ref{FIO})) We then put $\tilde{W}(u) = e^{i u P}$ and consider the conjugation $ \tilde{W}(u) \psi_V 1_{\Omega_o} \tilde{W}_{-1}(u).$  In a well-known way (cf. \cite{T}, \S 7.8), the complete  symbol expansion  of this operator is obtained recursively  by solving transport equations along the orbits of the Hamilton flow $\phi^u$ of $\sigma_P = H$.  The principal symbol   equals $\phi^{u *} ( \psi_V 1_{\Omega_o})$. Similarly, the   complete symbol at $\phi^1(x, \xi)$ of  $ \tilde{W} \psi_V 1_{\Omega_o} \tilde{W}_{-1}$   depends only on the germ of the complete symbol of  $\psi_V 1_{\Omega_o}$ at  $\phi^1(x, \xi)$ and on the germ of $\phi^1$ at $(x, \xi).$  Since $\phi^1(V_o) \subset V \cap T^*(int \Omega_o)$, the germ of the  complete symbol of $\psi_V 1_{\Omega_o}$ is identically one at $(x, \xi) \in V_o$ and hence the germ of the complete symbol of  $ \tilde{W} \psi_V 1_{\Omega_o} \tilde{W}_{-1}$ equals that of $\tilde{W} \tilde{W}_{-1}$ at $\phi(x, \xi)$,   hence equals one identically.  Since we have composed with $A$ satisfying $WF'(A) \subset V_o$ this calculation
is valid on microsupport of $A$ and thus $ \tilde{W} \psi_V  1_{\Omega_o} \tilde{W}_{-1} A \sim A$.

 The statement  $ \tilde{W} \; r \; \psi_V 1_{\Omega_o} \tilde{W}_{-1} A \sim 0$
follows from the fact that the underlying canonical transformation of this operator is $(\phi^{1})^{-1} r^* \phi^1$,
and the  graph of this canonical transformation  is disjoint from the 
$T^*(int \Omega_o) \cap V_o  \times T^*(int \Omega_o) \cap V_o.$  This follows since the graph of $r^*$ does not
intersect $T^*(int \Omega_o) \times T^*(int \Omega_o).$

Thus ${\cal E}(t)$ is a microlocal solution of the mixed Cauchy problem modulo
errors in $Op S^{2, K}$. 
 Since the Cauchy problem is well-posed, we have 
by Duhamel's principle that 
\begin{equation} ({\cal E}(t) - E(t) \sim  
\int_0^t  G_o(t - u) R(u) {\cal E} (u) du \end{equation}
where $R(u)$ is an error in the stated class and where 
where $G_o$ is the kernel of $\frac{ \sin t (\sqrt{\Delta_{\Omega}} )}{\sqrt{\Delta_{\Omega}}}$
i.e. of the mixed problem 
\begin{equation}\label{SIN} \begin{array}{ll} (\partial_t^2 - \bDelta) G_o = 0 & \mbox{on} \; 
\Omega_o \times \Omega_o \\ &  \\ 
G_o(t, \s, \y, \s', \y')|_{t = 0} = 0 & \frac{\partial }{\partial t} G_o(t, \s, \y, \s', \y')|_{t = 0} = Id  \\
& \\ G_o(t, \s, \y, \s', \y') = 0 & \mbox{for}\; \s = 0, L . \end{array} \end{equation}

Since we are working microlocally near $\gamma,$ we may replace $G_o$ modulo smoothing
operators by a Fourier integral parametrix (cf \cite{GM} \cite{PS}). Thus $\int_0^t  G_o(t - u) R(u)  {\cal E}_{\epsilon} Q )(u, \s, \y, \s', \y') du$ is a
sum of terms  of the form $A R_1 B$ with $R_1$ as above.
 \qed

\section{Wave invariants and normal form}

We now prove the crucial lemma that  the wave invariants of the normal form at $\gamma_o^m$, i.e. the  singularities of $Tr \cos t( F(|D|, \hat{I})$ at
$\s = 2 m L$, agree with the wave invariants  of $\sqrt{\Delta}$ at $\gamma^m$, i.e. the
 singularity of $Tr e^{i t \sqrt{\Delta}}$ at $t =2m L.$  Here, the $|D|$ in
 $F(|D|, \hat{I})$ refers to the Dirichlet $D$ on $[0, L]$ and the $\hat{I}$ as above refers to the homogeneous
action operator.

\begin{lem} Let $\gamma$ denote a non-degenerate bouncing ball orbit of $\Omega$ and let $\gamma_o$ denote
the corresponding bouncing ball orbit $\gamma_o$ of the normal form glow  of $\Omega_o$.  Then for all $m$, we have  $a_{\gamma^m k} (\sqrt{\Delta}) = a_{\gamma_o^m k} (F(|D|, \hat{I}).$ \end{lem}

\noindent{\bf Proof}

To prevent confusion between traces on $[0, L] \times \R$ and traces on $S^1 \times \R$ we will exclusively use the
notation 'res' for the latter space and will explicitly put in the cutoff $1_{\Omega_o}$ (the characteristic function
of $\Omega_o$) to indicate traces on  the former.  We have:

\begin{equation} \label{SONE}\begin{array}{l} a_{\gamma k} (\sqrt{\Delta_{\Omega}}) =
a_{\gamma_o k} (\sqrt{\bDelta_{\Omega_o}})  \\  \\
 = res \;  (D_t^k E(t)|_{t = 2 m L}) 1_{\Omega_o}.  \end{array} \end{equation}
The residue is given by  the integral of a density over
the fixed point set of the underlying  canonical relation. The
density is calculated using the  the method of stationary
phase on a manifold with boundary (see Theorem (\ref{POISSON}) or \cite{GM}),  essentially as in the boundaryless case. 
Thus the residue depends only on the kernel of $E(t)$ microlocally in the cone $V_o$ of
proposition (\ref{VZERO}).  Hence by lemma (\ref{PARA}) we have 
\begin{equation} \label{STWO}\begin{array}{l}res \; (D_t^k E(t)|_{t = 2 m L}) 1_{\Omega_o} =  2 res \; (D_t^k {\cal E}(t)|_{t = 2 m L}) 1_{\Omega_o} + 2 res\; (D_t^k R(t) |_{t = 2 m L}) 1_{\Omega_o} \end{array} \end{equation}
with $R(t) \in A_1 Op S^{2, K} A_2)$.  By Proposition (\ref{FIODROP}) we can drop the second term.  Then by the
 definition of ${\cal E}(t)$ we get that \eqref{STWO} equals 
\begin{equation}\label{STHREE} \begin{array}{l} 2 res 1_{\Omega_o} \tilde{W} (D_t^k   F_{o }(t)|_{t =  2mL}) \psi_V 1_{\Omega_o} \tilde{W}_{-1}.
  \end{array} \end{equation}

To
evalute this residue let us approximate $1_{\Omega_o}$ by a smooth cutoff $\tau_{\epsilon}(\s)$, equal
to one on $\Omega_o$ and supported for $\s \in (-\epsilon, \pi + \epsilon).$  
We have:
\begin{equation}\begin{array}{l} a_{\gamma^m k} (\sqrt{\Delta}) =  2 res \; 1_{\Omega_o} \tilde{W} (D_t^k   F_{o }(t)|_{t = 2m L} ) \psi_V 1_{\Omega_o} \tilde{W}_{-1}\\ \\ =2 \lim_{\epsilon \to 0} res\; \tau_{\epsilon}(x) \tilde{W} (D_t^k   F_{o }(t)|_{t = 2m L}) \psi_V 1_{\Omega_o} \tilde{W}_{-1}\\ \\
=2 \lim_{\epsilon \to 0} res\;  1_{\Omega_o}  \tilde{W}_{-1}  \tau_{\epsilon}  \tilde{W} (D_t^k  F_{o  }(t)  |_{t =  2m L}) \psi_V \end{array} \end{equation}
where in the last line we used the tracial property of $res$. 
We note  that 
$$  \tilde{W}_{-1} \tau_{\epsilon} \tilde{W} (D_t^k   F_{o }(t)  |_{t = 2m L}) \psi_V$$ is a (standard) Fourier integral operator on the boundaryless manifold  $S^1_{2L} \times \R$. 
As mentioned above, the  residue of $1_{\Omega_o}$ times this operator equals the integral of a density over
the fixed point set of the underlying  canonical relation. 
Since $  \tilde{W}_{-1} \tau_{\epsilon} \tilde{W}$
is a pseudodifferential operator, the canonical relation is simply that of $(D_t^k   F_{o }(t)  |_{t = 2m L})$.  
Due to the factor of  $1_{\Omega_o}$ the residue density may be calculated as if $\tau_{\epsilon} \equiv 1$. But then the factors of $\tilde{W}$ and $\tilde{W}_{-1}$ cancel and 
we are left with
 \begin{equation}\begin{array}{l} a_{\gamma^m k} (\sqrt{\Delta}) = 2 res\; 
1_{\Omega_o} (D_t^k   F_{o }(t)  |_{t = 2m L}) \psi_V . \end{array} \end{equation}
To complete the proof we must show that this residue equals 
\begin{equation} \label{CLAIM} res\; 
1_{\Omega_o} (D_t^k   F_{o }(t)  |_{t = 2m L}) \psi_V. \end{equation} The difference is just that the odd kernel equals $\Pi_o = \half (I -r)$ composed with the 
full  kernel.  The two terms of $\Pi_o$ give rise to 
 two components to the canonical relation of $F_{o}$.  They consist of  the graph of the Hamilton flow of the 
normal form , which we denote by  $G^{2 m L}_{nor}$,  and the graph of $G^{2mL}_{nor} \circ r.$  The latter map has no fixed points in $T*(int \Omega_o)$.  Hence  the residue of the
second term equals zero.  In the first term the factors of 2 cancel we give the result
claimed in (\ref{CLAIM}). 
 \qed

\begin{cor} The quantum Birkhoff normal coefficients, i.e. the coefficients of the
polynomials $p_k(\hat{I})$ are spectral invariants of $\Delta$. \end{cor}

The proof of the corollary from the lemma is identical to that in \cite{G} (see also \cite{Z.2}), so we omit the proof.

\subsection{Conclusion of Proof of Theorem}

We now complete the proof of the main result.  Since the quantum Birkhoff normal
coefficients of $\sqrt{\Delta}$ at $\gamma$ are spectral invariants, it follows
afortiori (by taking the interior symbol) that
the  Birkhoff normal form
for the metric Hamiltonian $H = |\xi|_g^2$ of the metric $g_{\Phi_0}$ at the bouncing ball
orbit is a spectral invariant.  Thus we may write:
\begin{equation} H = |\sigma| + \frac{\alpha}{L} I + \frac{p_1(I)}{|\sigma|} + \dots +
\frac{p_k(I)}{|\sigma|^k} + \dots \end{equation}
where $p_k$ is the homogeneous polynomial part of degree $k + 1$ in the quantum normal
form.  We now observe that this normal form of the Hamiltonian induces the Birkhoff normal
form of the Poincare map ${\cal P}_{\gamma}.$ 

Indeed, we first note that the coordinate $s$ is dual to $\sigma$ and hence is an angle
variable. Hence
the Hamilton flow in action-angle variables takes the form:
\begin{equation} \phi^t (s, \phi, |\sigma|, I) = (s + t \omega_{|\sigma|}, \phi + t
\omega_I, |\sigma|, I),\;\;\;\;\;\omega_{|\sigma|} = \partial_{|\sigma|} H, \omega_I = 
\partial_I H. \end{equation}
The billiard map from the bottom component of the boundary to the top component
is given by: 
\begin{equation} \beta_1 (\phi, I) = (\phi + t_1(\phi, |\sigma|, I) \omega_I, I) \end{equation}
where $t_1(\phi, |\sigma|, I)$ is the time until the trajectory through the initial vector
defined by  $(0, \phi,|\sigma|, I)$ hits the
upper part of the boundary.  We obviously have:
\begin{equation} t_1 (\phi, |\sigma|, I) = L/ \omega_{|\sigma|}. \end{equation}
Similarly, the billiard map for the return trip is given by
\begin{equation}
\beta_2(\phi, I) = (\phi + t(\phi, |\sigma|, I) \omega_I, I). \end{equation}
It follows that the Poincare map has the Birhoff normal form
\begin{equation} {\cal P}_{\gamma}(\phi, I) = (\phi + t_1 \omega_I + t_2 \omega_I \circ \beta_1,
I). \end{equation}

Thus, the Birkhoff normal form of ${\cal P}_{\gamma}$ is a spectral invariant. We now want
to apply the argument of Colin de Verdiere (\cite{CV}, \S 4,  THEOREME) to conclude that the domain $\Omega$ is determined by its spectrum.  Although Colin de Verdiere makes the assumption that $\gamma$ is a non-degenerate elliptic orbit, his argument is equally valid for hyperbolic
bouncing ball orbits. Let us briefly explain the necessary modifications.

Assuming that the non-degenerate bouncing ball orbit is the vertical axis, we write the graph of $\Omega$ over the horizontal axis as $y = f(x) = 1 + a_0 x^2 + \dots + a_n x^{2n - 2}.$ Only
even terms appear due to the left/right symmetry assumption. By definition, the Birkhoff normal form of the  (non-linear)
Poincare map ${\cal P}_{\gamma}$ is an expression for this map in local  action-angle variables $(I, \theta)$ on the transversal.  In the elliptic case we need to assume 
 (as above) that $\alpha/\pi \notin
\Q$ to ensure that the normal form exists;  no assumption is needed in the hyperbolic case. 
Following \cite{CV} we write  the normal form as
$$ T(I, \theta) = (I + O(I^{\infty}), \theta + b_0 + b_1 I + \cdots + O(I^{\infty})). $$
The key assertion is that there is an `upper-triangular'  bijection between the Taylor 
coefficients $a_j$ and the Birkhoff normal form  coefficients $b_k$.  This means (i) that
 $b_n = B + C a_n$
where $B,C$ depend only on $a_0, \dots, a_{n-1}$  and (ii) that $C \not= 0.$
Granted (i) - (ii),  $\{a_0, \dots, a_n\} \to
\{b_0, \dots, b_n\}$ can be inverted and the Taylor coefficients are determined by the normal
form coefficients.  

The proof of (i)--(ii) is almost exactly the same in the non-degenerate elliptic and
hyperbolic cases.  To convince the reader of this, we briefly recall the argument in
\cite{CV} and extend it to the hyperbolic case. 
The starting point is that  one may write down a
generating function $\phi$ for $T$ in terms of $f$:
$$ \mbox{graph}\; T = \{(x, \frac{\partial \phi}{\partial x};  x_1, - \frac{\partial \phi}{\partial x_1}: \frac{\partial \phi}{\partial s} = 0\}$$
with $\phi(x, x_1, s) = [(x -s)^2 + (f(s)^2]^{\half} + [(x_1 -s)^2 + (f(s)^2]^{\half}.$
Here, $x,x_1, s$ denote points on the horizontal axis.  In either the elliptic or hyperbolic
case, one may expand $\phi =2 +  \phi_0 + \cdots + \phi_n + \cdots$ where $\phi_0$ (resp. $\phi_n$) is  homogeneous of degree 2 (resp. degree $n + 2$) in $(x, x_1, s)$.  One easily
finds that $\phi_0(x, x_1, s) = \half (x^2 + x_1^2) - s(x + x_1) + (2 a_0 + 1) s^2)$ and
that $\phi_n = \psi_n + 2 a_n s^{2n+2}$ where $\psi_n$ depends only on $a_0, \dots, a_{n-1}.$
The linear Poincare map $P_{\gamma}$ has generating function $\phi_0$ and one finds that
$$P_{\gamma} = \left( \begin{array}{ll} A - 1 & -A \\ 2 - A & A - 1 \end{array} \right)$$
where $A = 2 (2 a_0 + 1).$  In the elliptic case $a_0 \in (-\half, 0)$ while in the hyperbolic
case $a_0 > 0$ or $a_0 < -\half.$  We assume $A \not= 0,$ i.e. $a_0 \not= -\half.$

Using the generating function one shows (\cite{CV}, Lemma 1) that 
\begin{equation} \label{CNF} T \left[ \begin{array}{l} x \\  \xi \end{array} \right] = T^0 \left[ \begin{array}{l} x \\ \xi \end{array} \right] + C (x - \xi)^{2n + 1} a_n \left[ \begin{array}{l} 1 \\  -1 \end{array} \right] + O(|(x, \xi)|^{2n + 2}) \end{equation}
where $T^0$ is the $(2n + 1)$ jet of $T$ computed with the assumption that $a_n = 0$ and
where $C \not= 0.$ The proof is a formal manipulation with Taylor series and only in the evaluation of $C$ does it matter whether the linear part of $T$ is elliptic or hyperbolic.
The constant $C$  has the form $C_1 A^{2n+2}$ where $C_1$ is universal.  Hence it is non-vanishing as long as $A \not= 0.$  

From the construction of the Birkhoff invariants (cf. \cite{SM}) and from (\ref{CNF})
it follows that $b_n = B + C a_n$.  It remains to show that $C \not= 0.$  To prove this,
we recall that two germs of  area-preserving transformations of $(\R^2, 0)$ for which
the eigenvalues of the linear parts are not roots of unity are symplectically equivalent
if and only if they have identical normal forms (\cite{SM}, p. 162).  To bring this down
to a finite dimensional statement, we define $G_{2n+1}$ as the group of (2n+1)-jets of
  area-preserving transformations of $(\R^2, 0)$, and $\Gamma_{2n+1} \subset G_{2n+1}$ as
the subgroup of elements of the form $Id + O(|(x,\xi)|^{2n+1}).$ Also let ${\cal O}_T =\{T' \in G_{2n+1}:
T' - T = O(|(x,\xi)|^{2n+1})\}$. Assuming the eigenvalues of the linear part of $T' \in {\cal O}_T$ are not
roots of unity, the orbit $\Gamma_{2n+1} \cdot T' \subset {\cal O}_T$ consists of elements $T'' \in {\cal O}_T$ with
the same Birkhoff normal form invariants up to and including $b_n$.  Put $T' = T(a_0, \dots, 0).$ 
If $C = 0$,  then $T(a_0, a_1, \dots, a_n)$ would on the orbit $\Gamma_{2n+1} \cdot T(a_0, \dots, 0)$  for
all $a_n$, i.e. the curve $a_n \to T(a_0, a_1, \dots, a_n)$ would be tangent to this orbit.
However, it is shown in \cite{CV} that it is transversal to the orbit.  In both elliptic and hyperbolic cases, 
we can choose coordinates $(u,v)$ in which the linear part of $T$ is diagonal and observe that 
the tangent vector to $a_n \to T(a_0, a_1, \dots, a_n)$  contains the monomials
$u^n v^{n+1}$ and $u^{n+1} v^n$ with non-zero coefficients.  However, these monomials cannot
occur in tangent vectors to $\Gamma_{2n+1} \cdot T(a_0, \dots, a_n)$. This statement
only involves the conjugation of the linear part of $T$ by $\Gamma_{2n+1}$ and is valid (with the
same proof) in both elliptic and hyperbolic cases. Therefore (i) - (ii) are valid in both cases and   $\{a_0, \dots, a_n\} \to \{b_0, \dots, b_n\}$ can be inverted. 
\qed


\begin{thebibliography}{HHHH}

\bibitem[A.M]{A.M} K.Andersson and R.B.Melrose, The propagation of singularities
along gliding rays, Invent.Math. 41 (1977), 23-95.

\bibitem[B.B]{B.B} V.M.Babic, V.S. Buldyrev: {\it Short-Wavelength Diffraction Theory},
Springer Series on Wave Phenomena 4, Springer-Verlag, New York (1991).

\bibitem[BM]{BM} L. Boutet de Monvel, Hypoelliptic operators with double characteristics
and related pseudodifferential operators, Comm. Pure Appl. Math. XXVII (1974), 585-639.

\bibitem[BMGH]{BMGH} L.Boutet de Monvel, A.Grigis, and B.Helffer, Parametrixes d'operateurs
pseudo-differentiels a caracteristiques multiples, Asterisque 34-35 (1976), 93-121.

\bibitem [CV]{CV} Y.Colin de Verdiere, Sur les longuers des trajectoires
periodiques d'un billard, In: P.Dazord and N. Desolneux-Moulis (eds.) {\it
Geometrie Symplectique et de Contact: Autour du Theoreme de Poincare-Birkhoff.
Travaux en Cours, Sem. Sud-Rhodanien de Geometrie III} Pairs: Herman (1984),
122-139.

\bibitem[DG]{DG} H. Duistermaat, V. Guillemin:  The spectrum of elliptic
operators and periodic bicharacteristics, Inv.\ Math.\ {\bf 29} (1975),
39--79.



\bibitem [G]{G} V.Guillemin, Wave trace invariants, Duke Math J. 83 (1996), 287-
352.

\bibitem[G.1]{G.1} V. Guillemin,  Wave-trace invariants and a theorem of Zelditch,
Duke Int.Math.Res.Not. 12 (1993), 303-308.



\bibitem[GM]{GM} V.Guillemin and R.B.Melrose, The Poisson summation formula for
manifolds with boundary, Adv.in Math. 32 (1979), 204 - 232.


\bibitem[H]{H} L. H\o"rmander, {\it The Analysis of Linear Partial Differential I  (III)
Operators}, Grund.Math.Wiss. 256 (275), Springer-Verlag, New York (1985).

\bibitem[L]{L} V.F.Lazutkin, Construction of an asymptotic series of
eigenfunctions of the ``bouncing ball'' type, Proc.Steklov Inst.Math. 95 (1968),
125- 140.

\bibitem[LT]{LT} V.F. Lazutkin and D.Ya.Terman, Number of quasimodes of `bouncing
ball' type, J.Soviet Math. (1984), 373-379.




\bibitem[PS]{PS} V.M.Petkov and L.N.Stoyanov, {\it Geometry of Reflecting Rays
and Inverse Spectral Problems}, John Wiley and Sons, N.Y. (1992).


\bibitem[SM]{SM} C.L.Siegel and J.Moser, {\it Lectures on Celestial Mechanics}, Grund.
Math.Wiss.Einz. 187, Springer-Verlag, New York (1971).

\bibitem[T]{T} M.E. Taylor, {\it Partial Differential Equations II}, Appl. Math.Sci. 116, Springer-Verlag (1996).



\bibitem[Z.1]{Z.1} S.Zelditch,  Wave invariants at elliptic closed geodesics,
 GAFA 7 (1997), 145-213.

\bibitem[Z.2]{Z.2} ----------,  Wave invariants for non-degenerate closed geodesics,
GAFA 8 (1998), 179-217.

\bibitem[Z.3]{Z.3} ----------, The inverse spectral problem for surfaces of
revolution, J. Diff. Geom. 49 (1998), 207-264. 

\bibitem[Z.4]{Z.4} ----------, Normal form of the wave group and inverse spectral theory, {\it Journees Equations
aux Derivees Partielles}, Saint-Jean-de-Monts (1998).

\bibitem[Z.5]{Z.5} ----------, Spectral determination of analytic bi-axisymmetric plane
domains (announcement), (preprint 1998).



\end{thebibliography}
\end{document}